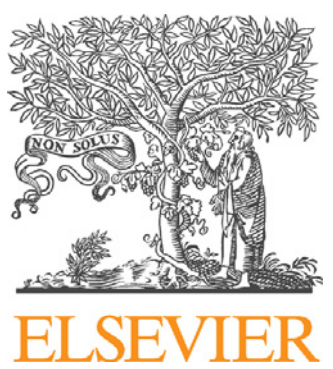

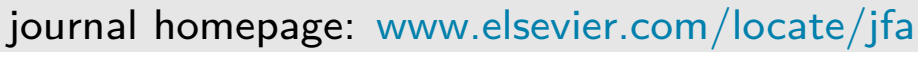

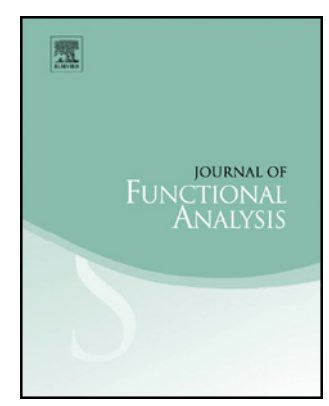

Regular Article

# On the lack of selection for the transport equation over a dense set of vector fields

Jules Pitcho [a,b,*]

[a] *ENS de Lyon, UMPA, 46 allée d'Italie, 69364 Lyon, France*
[b] *GSSI, Via Michele Iacobucci, 2, 67100 L'Aquila, Italy*

ABSTRACT

We construct a set of bounded vector fields dense in $L^p((0,2); W^{s,p}_{loc}(\mathbb{R}^d;\mathbb{R}^d))$ for $1 \leq p < +\infty$ and $0 \leq s < 1$ with $p < 1/s$ for which smooth regularisation of the vector field does not give a selection criterion for the continuity equation, thereby showing that the two examples constructed in [10,12] are generic.

## 1. Introduction

We study the initial value problem for the continuity equation

$$\begin{cases} \partial_t \rho + \operatorname{div}_x(\boldsymbol{b}\rho) = 0, \\ \rho(0,x) = \bar{\rho}(x), \end{cases} \tag{IVP}$$

* Correspondence to: GSSI, Via Michele Iacobucci, 2, 67100 L'Aquila, Italy.
*E-mail address:* jules.pitcho@gssi.it.

posed on $[0, T] \times \mathbb{R}^d$, where $\boldsymbol{b} = \boldsymbol{b}(t, x)$ is a given vector field, $\rho = \rho(t, x)$ is an unknown real-valued function, and $\mathrm{div}_x$ is the divergence operator on vector fields on $\mathbb{R}^d$. We are interested in bounded weak solutions of (IVP).

**Definition 1.1.** Consider a bounded vector field $\boldsymbol{b} : [0, T] \times \mathbb{R}^d \to \mathbb{R}^d$ and an initial datum $\bar{\rho} \in L^\infty(\mathbb{R}^d)$. We shall say that $\rho \in C([0, T]; w^* - L^\infty(\mathbb{R}^d))$ is a bounded weak solution to (IVP) along $\boldsymbol{b}$, if for every $\phi \in C_c^\infty((0, T) \times \mathbb{R}^d)$

$$\int_0^{+\infty} \int_{\mathbb{R}^d} \rho \left( \frac{\partial \phi}{\partial t} + \boldsymbol{b} \cdot \nabla_x \phi \right) \, dx dt = 0, \qquad \text{and} \qquad \rho(0, \cdot) = \bar{\rho}(\cdot).$$

Existence of a bounded weak solution holds, if we assume that $[\mathrm{div}_x \boldsymbol{b}]^- \in L^1((0, T); L^\infty(\mathbb{R}^d))$ by a regularisation and compactness argument. When the vector field $\boldsymbol{b}$ is sufficiently rough, uniqueness of bounded weak solutions fails, even for divergence-free vector fields, as the example of Depauw [13] shows. The following question may be of interest: *is there a robust selection criterion for bounded weak solutions of* (IVP) *assuming that* $\boldsymbol{b}$ *is only bounded and divergence-free?* When $\boldsymbol{b}$ is in $L^1_t BV_x$, the work of Ambrosio [1] following on that of DiPerna and Lions [14] answers this question affirmatively. In fact, uniqueness of bounded weak solution of (IVP) was proved in [1], assuming that $\boldsymbol{b}$ is in $L^1_t BV_x$, and that $\mathrm{div}_x \boldsymbol{b}$ is absolutely continuous with respect to the Lebesgue measure for almost every time.

As (IVP) is well-posed when $\boldsymbol{b}$ is smooth and bounded, a first naive attempt at answering the above question is to regularise the vector field $\boldsymbol{b}$, and to study limit points of solutions of the regularised equation, in the hope that, through uniqueness of the limit point, a single solution of (IVP) will be selected. For a specific subclass of bounded, divergence-free vector fields, which includes vector fields for which uniqueness of bounded weak solutions of (IVP) fails, such a criterion was proven by the author in [20] using regularisations by convolution of $\boldsymbol{b}$. He also presents a Lagrangian counterpart to this selection result in [19], where interestingly he shows that for Depauw's example [13], a probability measure concentrated on several distinct integral curves is selected for Lebesgue almost every initial datum. When uniqueness of bounded weak solutions holds, as it does in Ambrosio's theorem, then regularisations of $\boldsymbol{b}$ generating solutions of (IVP) which are uniformly bounded in $L^\infty_{t,x}$ automatically converge to the unique bounded weak solution.

However, regularisation of the vector field is, in general, not a selection criterion. This was first observed by Ciampa, Crippa and Spirito in [10], where they constructed an unbounded, divergence-free vector field for which two distinct regularisations yield two limiting solutions. Afterwards, in [12] De Lellis and Giri constructed a bounded, divergence-free vector field, as well as an initial datum for which two distinct regularisations yield two distinct limiting solutions. Colombo, Crippa and Sorella then constructed in [11] vectors fields in $C^\alpha([0, 2] \times \mathbb{R}^2)$ where $\alpha < 1$ is fixed but arbitrary, and corre-sponding initial data for which regularisation by convolution of the vector field is not a selection criterion.

Another regularisation scheme consists in regularising the equation by adding diffusivity, namely for a fixed value of the parameter $\nu \in (0,1)$, one considers the following initial value problem

$$\begin{cases} \partial_t \rho^\nu + \operatorname{div}_x(\boldsymbol{b}\rho^\nu) - \nu \Delta_x \rho^\nu = 0, \\ \rho^\nu(0,x) = \bar{\rho}(x), \end{cases} \tag{$\nu$-IVP}$$

which is well-posed for bounded initial data. We then say that *the diffusivity regularisation scheme gives a selection criterion*, if the family $\{\rho^\nu\}_{\nu\in(0,1)}$ of unique solutions of ($\nu$-IVP) admits a unique weak limit point as $\nu \downarrow 0$. In [11], the authors constructed bounded, divergence-free vector fields and corresponding bounded initial data for which the vanishing diffusivity regularisation scheme does not give a selection criterion. In [16], the authors further constructed a bounded, divergence-free vector field for which, for every bounded initial data, the diffusivity regularisation scheme is not a selection criterion. For a class of divergence-free vector fields for which uniqueness of bounded weak solutions fails, it is shown in [18] that the diffusivity regularisation scheme gives a selection criterion for every bounded initial datum.

A question related to selection by the diffusivity regularisation scheme is anomalous dissipation. We say that a vector field $\boldsymbol{b}$ and a bounded initial datum $\bar{\rho}$ exhibit anomalous dissipation, if the family $\{\rho^\nu\}_{\nu\in(0,1)}$ of unique solutions of ($\nu$-IVP) satisfies for some $R > 0$ that

$$\limsup_{\nu\downarrow 0} \nu \int_0^T \int_{[-R,R]^d} |\nabla \rho^\nu(t,x)|^2 dxdt > 0.$$

The vector fields and corresponding initial data constructed in [11] for which the vanishing diffusivity scheme is not a selection criterion are shown to exhibit anomalous dissipation. In [21], it is shown that anomalous dissipation implies non-uniqueness of the backward advection equation. But for vector fields which have anomalous dissipation, there may be a unique vanishing diffusivity solution as the construction of shows [15]. See also [4,9,17] for works on anomalous dissipation.

We here argue that the lack of selection demonstrated by De Lellis and Giri is generic, in analogy with the nonlinear wave equation, where Sun and Tzvetkov [22] showed ill-posedness on a dense set of initial data of super-critical regularity. We will thus construct a dense set of bounded vector fields for which the lack of selection demonstrated by De Lellis and Giri holds. In addition to being dense, this set will have the following property, which expresses that there is no distinguished time, when measuring the vector field in the $L^\infty_x$ norm.

**Definition 1.2.** We shall say that a subset $D$ of $L^\infty([0,T];L^\infty(\mathbb{R}^2;\mathbb{R}^2))$ is *not bounded at any time in* $[0,T]$, if for every $\gamma > 0$ there exists $\boldsymbol{b} \in D$ such that for every interval $I \subset [0,T]$ with non-empty interior, we have $\|\boldsymbol{b}\|_{L^\infty(I;L^\infty_x)} \geq \gamma$.

The above definition forces the set $D$ to include vector fields which are arbitrarily large over time interval. Without this requirement, the main theorem of this paper could be proved by another construction perturbing the basic building blocks described in Section 3 away from the singularity. At any rate, we will refer to a smooth and bounded sequence $(\boldsymbol{b}^q)_{q\in\mathbb{N}}$ such that $\boldsymbol{b}^q \to \boldsymbol{b}$ in $L^1_{loc}$ as $q \to +\infty$ as a regularisation of $\boldsymbol{b}$, and we now state our main theorem.

**Theorem 1.3.** *Let* $1 \leq p < +\infty$ *and* $0 \leq s < 1$ *such that* $p < 1/s$. *There exists a subset* $D \subset L^\infty([0,2];L^\infty(\mathbb{R}^2;\mathbb{R}^2))$, *which is dense in* $L^p((0,2);W^{s,p}_{loc}(\mathbb{R}^2;\mathbb{R}^2))$, *and which is not bounded at any time in* $[0,2]$ *such that the following holds. For every* $\boldsymbol{b} \in D$, *there exists an initial datum* $\bar{\rho} \in L^\infty(\mathbb{R}^2)$, *two bounded weak solutions* $\rho$ *and* $\tilde{\rho}$ *of* (IVP) *along* $\boldsymbol{b}$ *with initial datum* $\bar{\rho}$, *and two regularisations* $(\boldsymbol{b}^q)_{q\in\mathbb{N}}$ *and* $(\tilde{\boldsymbol{b}}^q)_{q\in\mathbb{N}}$ *of* $\boldsymbol{b}$ *such that:*

- *the unique solutions of* (IVP) *along* $\boldsymbol{b}^q$ *with initial datum* $\bar{\rho}$ *converge in* $C([0,2];w^*-L^\infty(\mathbb{R}^2))$ *to* $\rho$ *as* $q \to +\infty$*;*
- *the unique solutions of* (IVP) *along* $\tilde{\boldsymbol{b}}^q$ *with initial datum* $\bar{\rho}$ *converge in* $C([0,2];w^*-L^\infty(\mathbb{R}^2))$ *to* $\tilde{\rho}$ *as* $q \to +\infty$.

We stress that the initial datum $\bar{\rho}$ depends on the vector field $\boldsymbol{b}$ in $D$. We also note that Bianchini and Zizza in [7] studied the density properties of divergence-free vector fields $\boldsymbol{b}$ in $L^1((0,1);BV(\mathbb{T}^2;\mathbb{R}^2))$, whose flow is mixing in the sense that the compositions of the flow map from time 0 to time 1 satisfy a mixing property in the sense of dynamical systems.

**Remark 1.4.** It would be interesting to construct a dense set of vector fields for which the diffusivity regularisation scheme does not give a selection criterion, and for which regularisation by convolution does not give a selection criterion. It would also be interesting to construct a dense set of autonomous vector fields in spatial dimension $d \geq 3$ for which there is no selection by regularisation as in the above theorem. Finally, it would also be interesting to investigate whether such sets of vector fields can be $G_\delta$ dense.

### *1.1. Outline of ideas*

The existence of vector fields $(\boldsymbol{b}_\lambda)_{\lambda\in\mathbb{N}}$, which converge to zero as $\lambda \to +\infty$, and corresponding initial data for which two different regularisations select two distinct bounded weak solutions of (IVP) follow from a rescaling of the building blocks of [12]. The important property of the building blocks is that they have a singularity at time 1, where the $BV$ norm of the vector field is not integrable, and where the solution becomes com-pletely mixed. From that time onwards, the solution can either remain completely mixed, or unmix leading to a non-uniqueness mechanism. Depending on how the vector field is regularised around this singular time will yield either the unmixing solution or the solution which stays completely mixed in the limit.

To prove density, we will add an arbitrary smooth perturbation $\boldsymbol{w}$, which need not be vanishing at any time, while simultaneously performing a transform on the vector field $\boldsymbol{b}_\lambda$ using the flow of $\boldsymbol{w}$ to obtain a family $(\boldsymbol{b}_{\lambda,\boldsymbol{w}})_{\lambda\in\mathbb{N}}$ (see Section 4.1). The vector fields $\boldsymbol{b}_{\lambda,\boldsymbol{w}}$ can have divergence with non-vanishing singular part with respect to the Lebesgue measure. The cost in $L^p_t W^{s,p}_x$ of this transform will be under control by smooth norms of $\boldsymbol{w}$, whereby the family $(\boldsymbol{b}_{\lambda,\boldsymbol{w}})_{\lambda\in\mathbb{N}}$ will converge to $\boldsymbol{w}$ in $L^p_t W^{s,p}_x$ as $\lambda \to +\infty$. The key point is choosing the transform on $\boldsymbol{b}_\lambda$ to decouple the flow at large scales and at small scales, so that the flow of the vector field $\boldsymbol{b}_{\lambda,\boldsymbol{w}}$ is given by the composition of the flow of the large scales with the flow of the small scales, and upon suitable truncation of $\boldsymbol{b}_{\lambda,\boldsymbol{w}}$ (see Section 4.2), the flow may be uniquely identified thanks to the uniqueness result of Bianchini and Bonicatto [5] for the flow of a *nearly incompressible* vector field in $L^1_t BV_x$ (see Section 2).

It would also be interesting to give a construction of a set $D$ of vector fields as in Theorem 1.3 via a different method.

### 1.2. Plan of the paper

In Section 2, we collect some useful results concerning the well-posedness of the continuity equation in the setting of $BV$ vector fields. We then collect the corresponding results on the flow of these vector fields. In Section 3, we first construct vector fields $(\boldsymbol{b}_\lambda)_{\lambda\in\mathbb{N}}$, which converge to zero as $\lambda \to +\infty$, and corresponding initial data for which two different truncations procedures select two distinct bounded weak solutions. In Section 4, we add a smooth perturbation $\boldsymbol{w}$ to the vector fields $(\boldsymbol{b}_\lambda)_{\lambda\in\mathbb{N}}$ and perform a transform on $\boldsymbol{b}_\lambda$ to obtain vector fields $(\boldsymbol{b}_{\lambda,\boldsymbol{w}})_{\lambda\in\mathbb{N}}$, which converge to $\boldsymbol{w}$ as $\lambda \to +\infty$. We then exhibit two different truncations of the singularity of $\boldsymbol{b}_{\lambda,\boldsymbol{w}}$ along which (IVP) has two distinct solutions, which in the limit where the truncation is removed, yield two distinct solutions of (IVP) along $\boldsymbol{b}_{\lambda,\boldsymbol{w}}$. In Section 5, we conclude the proof of Theorem 1.3.

### Acknowledgments

The author is thankful to his advisor Nikolay Tzvetkov for his support, for bringing the question of genericity of non-uniqueness constructions to his attention, and for useful comments. The author acknowledges the hospitality of the Pitcho Centre for Scientific Studies where part of this work was carried out. The author is thankful to the anonymous referees for numerous useful suggestions, which have improved this manuscript. The author has been partially funded by Simons Foundation Award ID: 651675.

## 2. Preliminaries

In this section $I$ is a closed interval whose interior we denote by $\mathring{I}$, and $\boldsymbol{b} : I \times \mathbb{R}^d \to \mathbb{R}^d$ is a Borel vector field whose essential supremum is bounded.

### 2.1. A boundary value problem

For an arbitrary $s \in I$, we consider the following boundary value problem posed on $I \times \mathbb{R}^d$:

$$\begin{cases} \partial_t \rho + \mathrm{div}_x(\boldsymbol{b}\rho) = 0, \\ \qquad\qquad \rho(s,x) = \bar{\rho}(x). \end{cases} \tag{BVP}$$

**Definition 2.1.** Consider a boundary datum $\bar{\rho} \in L^\infty(\mathbb{R}^d)$. We shall say $\rho \in C(I; w^* - L^\infty(\mathbb{R}^d))$ is a bounded weak solution to (BVP) along $\boldsymbol{b}$ with boundary datum $\bar{\rho}$ at time $s$, if for every $\phi \in C_c^\infty(\mathring{I} \times \mathbb{R}^d)$

$$\int_{\mathbb{R}} \int_{\mathbb{R}^d} \rho \left( \frac{\partial \phi}{\partial t} + \boldsymbol{b} \cdot \nabla_x \phi \right) dxdt = 0 \qquad \text{and} \qquad \rho(s,\cdot) = \bar{\rho}(\cdot).$$

Existence of bounded weak solutions of (BVP) follows from a regularisation and compactness argument, see for instance [2]. Nearly incompressible vector fields are a weaker version of vector fields with $\mathrm{div}_x\, \boldsymbol{b} \in L^1(I; L^\infty(\mathbb{R}^d))$. They are defined as follows.

**Definition 2.2.** We shall say that $\boldsymbol{b}$ is *nearly incompressible*, if there exists a constant $C > 0$ and a bounded weak solution $\rho$ of (IVP) along $\boldsymbol{b}$ such that $C^{-1} < \rho(t,x) < C$ for every $t \in I$ and a.e. $x \in \mathbb{R}^d$.

A regularisation and compactness argument shows that if $\mathrm{div}_x\, \boldsymbol{b} \in L^1(I; L^\infty(\mathbb{R}^d))$, then $\boldsymbol{b}$ is nearly incompressible. The converse is not true as the example in [6] shows. Notice that for $J$ a closed subinterval of $I$, if $\boldsymbol{b}$ is nearly incompressible, then the restriction of $\boldsymbol{b}$ to $J \times \mathbb{R}^d$ is also nearly incompressible. The following theorem can be extracted from the important work of Bianchini and Bonicatto [5] and concerns uniqueness and stability of bounded weak solutions of (BVP) along a nearly incompressible vector field.

**Theorem 2.3.** *Consider a boundary datum $\bar{\rho} \in L^\infty(\mathbb{R}^d)$. Assume that $\boldsymbol{b} \in L^1(I; BV_{loc}(\mathbb{R}^d; \mathbb{R}^d))$ and that $\boldsymbol{b}$ is nearly incompressible. Then, the following holds:*

- (i) *(uniqueness) there exists a unique bounded weak solution $\rho$ of (BVP) along $\boldsymbol{b}$, which belongs to $C(I; L^1_{loc}(\mathbb{R}^d))$;*
- (ii) *(stability) for every regularisation $(\boldsymbol{b}^q)_{q \in \mathbb{N}}$ of $\boldsymbol{b}$ such that there exists $C > 0$ such that the flow $\boldsymbol{X}^q$ along $\boldsymbol{b}^q$ satisfies*

$$C^{-1}\mathscr{L}^d \le \boldsymbol{X}^q(t,\cdot)_\#\mathscr{L}^d \le C\mathscr{L}^d \qquad \forall q \in \mathbb{N}, \tag{2.1}$$

*the unique bounded weak solution $\rho^q$ of (BVP) along $\boldsymbol{b}^q$ with boundary datum $\bar{\rho}$ converges to $\rho$ in $C(I; L^1_{loc}(\mathbb{R}^d))$ as $q \to +\infty$.*

Notice that the existence of a regularisation $(\boldsymbol{b}^q)_{q\in\mathbb{N}}$ whose flows satisfy (2.1) implies that $\boldsymbol{b}$ is nearly incompressible.

*2.2. Regular Lagrangian flows*

The regular Lagrangian flow provides a robust measure-theoretic notion of a flow of a rough vector field. For later use, we here define regular Lagrangian flows starting at time $s$.

**Definition 2.4.** Consider a real number $s \in I$. We shall say that a Borel map $\boldsymbol{X} : I\times\mathbb{R}^d \to \mathbb{R}^d$ is a regular Lagrangian flow along $\boldsymbol{b}$ starting at $s$, if

(i) for $\mathscr{L}^d$-a.e. $x \in \mathbb{R}^d$, $t \mapsto \boldsymbol{X}(t,x)$ is an absolutely continuous solution of $\partial_t \boldsymbol{X}(t,x) = \boldsymbol{b}(t,\boldsymbol{X}(t,x))$ with $\boldsymbol{X}(s,x) = x$, namely

$$\boldsymbol{X}(t,x) - x = \int_s^t \boldsymbol{b}(u,\boldsymbol{X}(u,x))du;$$

(ii) there exists a constant $C > 0$ called a compressibility constant such that for every $t \in I$, we have $C^{-1}\mathscr{L}^d \le \boldsymbol{X}(t,\cdot)_\#\mathscr{L}^d \le C\mathscr{L}^d$.

A Borel map $\boldsymbol{X} : I \times \mathbb{R}^d \to \mathbb{R}^d$ will be called *a unique regular Lagrangian flow* along $\boldsymbol{b}$ starting at $s$, if for every $\tilde{\boldsymbol{b}}$, which coincides $\mathscr{L}^{d+1}$-a.e. with $\boldsymbol{b}$, we have that for every regular Lagrangian flow $\tilde{\boldsymbol{X}}$ along $\tilde{\boldsymbol{b}}$ starting at $s$, it holds that $\boldsymbol{X}(\cdot,x) = \tilde{\boldsymbol{X}}(\cdot,x)$ for $\mathscr{L}^d$-a.e. $x \in \mathbb{R}^d$. By a slight abuse of language, a unique regular Lagrangian flow will also be referred to as *the* unique regular Lagrangian flow. The following theorem can be extracted from the work of Bianchini and Bonicatto [5], which we formulate as the counterpart, at the level of the flow, of the uniqueness in Theorem 2.3.

**Theorem 2.5.** *Consider a real number $s \in I$. Assume that $\boldsymbol{b} \in L^1(I; BV_{loc}(\mathbb{R}^d;\mathbb{R}^d))$, and that $\boldsymbol{b}$ is nearly incompressible. Then, there exists a unique regular Lagrangian flow along $\boldsymbol{b}$ starting at $s$.*

*2.3. Estimates on the flow of a smooth vector field*

Let $\boldsymbol{w} : I \times \mathbb{R}^d \to \mathbb{R}^d$ be a vector field in $L^\infty(I; C^1(\mathbb{R}^d;\mathbb{R}^d))$, let $s \in I$, and let $\boldsymbol{X}_{\boldsymbol{w}}$ be the unique flow along $\boldsymbol{w}$ starting at time $s$, namely $\boldsymbol{X}_{\boldsymbol{w}}$ solves

$$\begin{cases} \partial_t \boldsymbol{X}_{\boldsymbol{w}}(t,x) = \boldsymbol{w}(t, \boldsymbol{X}_{\boldsymbol{w}}(t,x)), \\ \quad \boldsymbol{X}_{\boldsymbol{w}}(s,x) = x, \end{cases} \tag{2.2}$$

in the sense of distributions over $I$. We define the Jacobian as $J\boldsymbol{X}_{\boldsymbol{w}}(t,x) = \det D_x \boldsymbol{X}_{\boldsymbol{w}}(t,x)$. A computation shows that

$$\begin{cases} \partial_t J\boldsymbol{X}_{\boldsymbol{w}}(t,x) = \operatorname{div}_x \boldsymbol{w}(t, \boldsymbol{X}_{\boldsymbol{w}}(t,x)) J\boldsymbol{X}_{\boldsymbol{w}}(t,x), \\ \quad J\boldsymbol{X}_{\boldsymbol{w}}(s,x) = 1. \end{cases} \tag{2.3}$$

Therefore, Grönwall inequality implies

$$\begin{aligned} \exp\Big[ -\Big| \int_s^t \| [\operatorname{div}_x \boldsymbol{w}(s,\cdot)]^- \|_{L^\infty_x} ds \Big| \Big] \mathscr{L}^d &\leq J\boldsymbol{X}_{\boldsymbol{w}}(t,x) \\ &\leq \exp\Big[ \Big| \int_s^t \| [\operatorname{div}_x \boldsymbol{w}(s,\cdot)]^+ \|_{L^\infty_x} ds \Big| \Big] \mathscr{L}^d. \end{aligned} \tag{2.4}$$

By a change of variables, we have

$$\boldsymbol{X}_{\boldsymbol{w}}(t,\cdot)_{\#} \mathscr{L}^d = \frac{1}{J\boldsymbol{X}_{\boldsymbol{w}}(t,\cdot)} \circ \boldsymbol{X}_{\boldsymbol{w}}^{-1}(t,\cdot) \mathscr{L}^d, \tag{2.5}$$

which implies, in view of (2.4), that

$$\begin{aligned} \exp\Big[ -\Big| \int_s^t \| [\operatorname{div}_x \boldsymbol{w}(s,\cdot)]^+ \|_{L^\infty_x} ds \Big| \Big] \mathscr{L}^d &\leq \boldsymbol{X}_{\boldsymbol{w}}(t,\cdot)_{\#} \mathscr{L}^d \\ &\leq \exp\Big[ \Big| \int_s^t \| [\operatorname{div}_x \boldsymbol{w}(s,\cdot)]^- \|_{L^\infty_x} ds \Big| \Big] \mathscr{L}^d. \end{aligned} \tag{2.6}$$

For a proof, we refer to [2]. We also have that for every $k \in \mathbb{N}$, there exists a constant $C_{k,d} > 0$ such that we have

$$\| \boldsymbol{X}_{\boldsymbol{w}}(t,\cdot) \|_{C^k_x} \leq C_{k,d} \exp\Big( \Big| \int_s^t \| \boldsymbol{w}(u,\cdot) \|_{C^k_x} du \Big| \Big). \tag{2.7}$$

## 3. Unperturbed vector fields

We describe the basic building blocks, which we shall use to construct a family of vector fields, which accumulates at the zero. We introduce the parameter $\lambda \in \mathbb{N}$. The building block $\boldsymbol{b}_\lambda$ will be at scale $2^{-\lambda}$. These building blocks, based on a mixing-unmixing mechanism are adapted from [12] to form a set $\{\boldsymbol{b}_\lambda : \lambda \in \mathbb{N}\}$, which accumulates at zero in suitable topologies. This section is based on [12] with the additional parameter $\lambda$.

### 3.1. The basic building blocks

We construct a family of bounded, divergence-free vector field $\boldsymbol{b}_\lambda : [0,2] \times \mathbb{R}^2 \to \mathbb{R}^2$, as well as two bounded weak solutions of $\zeta_\lambda$ and $\tilde{\zeta}_\lambda$ of (BVP) along $\boldsymbol{b}_\lambda$ with boundary datum $\bar{\rho} = 1/2$ at time 1.

Introduce the following two lattices on $\mathbb{R}^2$, namely $\mathcal{L}^1 := \mathbb{Z}^2 \subset \mathbb{R}^2$ and $\mathcal{L}^2 := \mathbb{Z}^2 + (\frac{1}{2}, \frac{1}{2}) \subset \mathbb{R}^2$. To each lattice, associate a subdivision of the plane into squares, which have vertices lying in the corresponding lattices, which we denote by $\mathcal{S}^1$ and $\mathcal{S}^2$. Then consider the rescaled lattices $\mathcal{L}^1_\lambda := 2^{-\lambda}\mathbb{Z}^2$ and $\mathcal{L}^2_\lambda := (2^{-\lambda-1}, 2^{-\lambda-1}) + 2^{-\lambda}\mathbb{Z}^2$ and the corresponding square subdivision of $\mathbb{R}^2$, respectively $\mathcal{S}^1_\lambda$ and $\mathcal{S}^2_\lambda$. The centres of the squares $\mathcal{S}^1_\lambda$ are elements of $\mathcal{L}^2_\lambda$ and vice-versa.

Next, define the following 2-dimensional autonomous vector field:

$$\boldsymbol{v}(x) = \begin{cases} (0, 4x_1)^t \ , \text{ if } 1/2 > |x_1| > |x_2| \\ (-4x_2, 0)^t \ , \text{ if } 1/2 > |x_2| > |x_1| \\ (0,0)^t \ , \text{ otherwise.} \end{cases}$$

$\boldsymbol{v}$ is a bounded, divergence-free vector field, whose derivative is a finite matrix-valued Radon measure given by

$$D\boldsymbol{v}(x_1,x_2) = \begin{pmatrix} 0 & 0 \\ 4\mathrm{sgn}(x_1) & 0 \end{pmatrix} \mathscr{L}^2 \llcorner_{\{|x_2|<|x_1|<1/2\}} + \begin{pmatrix} 0 & -4\mathrm{sgn}(x_2) \\ 0 & 0 \end{pmatrix} \mathscr{L}^2 \llcorner_{\{|x_1|<|x_2|<1/2\}}$$
$$+ \begin{pmatrix} 4x_2\mathrm{sgn}(x_1) & -4x_2\mathrm{sgn}(x_2) \\ 4x_1\mathrm{sgn}(x_1) & -4x_1\mathrm{sgn}(x_2) \end{pmatrix} \mathscr{H}^1 \llcorner_{\{|x_1|=|x_2|, 0<|x_1|,|x_2|\le 1/2\}}$$

Periodise $\boldsymbol{v}$ by defining $\Lambda_\lambda = \{(y_1,y_2) \in 2^{-\lambda}\mathbb{Z}^2 : 2^\lambda(y_1+y_2) \text{ is even}\}$ and setting

$$\boldsymbol{u}_\lambda(x) = 2^{-\lambda} \sum_{y \in \Lambda_\lambda} \boldsymbol{v}(2^\lambda x - y) \, .$$

Even though $\boldsymbol{u}_\lambda$ is non-smooth, it is in $BV_{loc}(\mathbb{R}^2;\mathbb{R}^2)$.

By Theorem 2.5, there exists a unique regular Lagrangian flow $\boldsymbol{X}_\lambda$ along $\boldsymbol{u}_\lambda$, which can be described explicitly.

(R) The map $\boldsymbol{X}_\lambda(t,\cdot)$ is Lipschitz on each square $S$ of $\mathcal{S}^2_\lambda$ and $\boldsymbol{X}_\lambda(1/2,\cdot)$ is a clockwise rotation of $\pi/2$ radians of the "filled" $S$, while it is the identity on the "empty ones". In particular for every $j \ge \lambda$, $\boldsymbol{X}_\lambda(1/2,\cdot)$ maps an element of $\mathcal{S}^1_j$ rigidly onto another element of $\mathcal{S}^1_j$. For $j = \lambda$ we can be more specific. Each $S \in \mathcal{S}^2_\lambda$ is formed precisely by 4 squares of $\mathcal{S}^1_\lambda$: in the case of "filled" $S$ the 4 squares are permuted in a 4-cycle clockwise, while in the case of "empty" $S$ the 4 squares are kept fixed (Fig. 1).

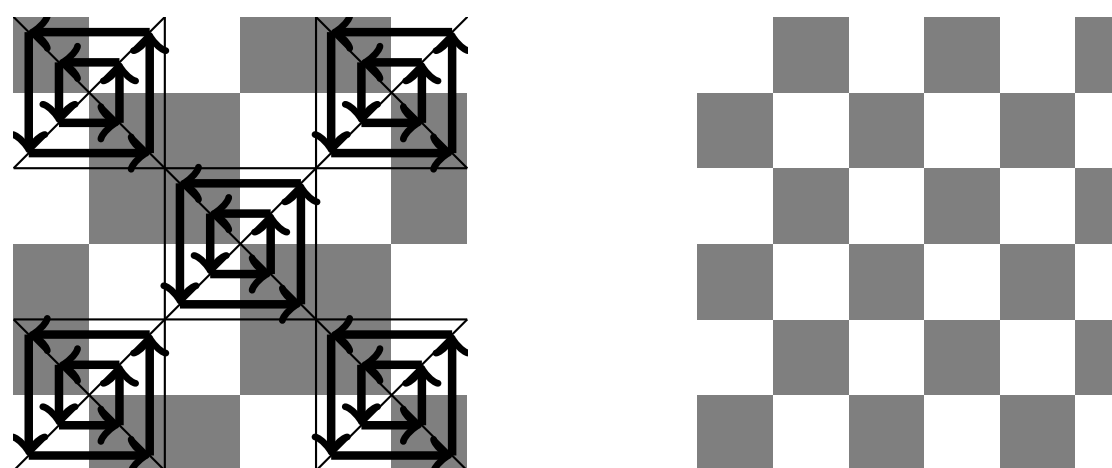

**Fig. 1.** Action of the flow of $\boldsymbol{u}_\lambda$ from $t = 0$ to $t = 1/2$. The shaded region denotes the set $\{\zeta_1 = 1\}$. Credits to [12] for this figure.

We define $\boldsymbol{b}_\lambda$ as follows. Set $\boldsymbol{b}_\lambda(t,x) = \boldsymbol{u}_\lambda(x)$ for $0 < t < 1/2$ and $\boldsymbol{b}_\lambda(t,x) = \boldsymbol{u}_\lambda(2^k x)$ for $1 - 1/2^k < t < 1 - 1/2^{k+1}$. Finally, set $\boldsymbol{b}_\lambda(t,x) = \boldsymbol{b}_\lambda(2-t,x)$. Notice that by construction we have

$$\|\boldsymbol{b}_\lambda\|_{L^\infty([0,2]\times\mathbb{R}^2;\mathbb{R}^2)} \leq 2^{-\lambda}\|\boldsymbol{v}\|_{L^\infty(\mathbb{R}^2;\mathbb{R}^2)}. \tag{3.1}$$

We now make the following observation due to the self-similarity of the contruction.

(O) Consider a bounded weak solution $\rho$ of (BVP) along $\boldsymbol{b}_\lambda$. As a direct consequence of the self-similarity of the construction, we have the following observation. If for some $q \in \mathbb{N}$ and some $\alpha \in \mathbb{R}$, we have $⨍_S \rho(1+2^{-q},x)dx = \alpha$ for every $S \in \mathcal{S}^1_{\lambda q}$, then for every $q' \leq q$, we have $⨍_S \rho(1+2^{-q'},x)dx = \alpha$ for every $S \in \mathcal{S}^1_{\lambda q}$.

Let $\bar{\zeta}_\lambda(x) = \lfloor 2^\lambda x_1 \rfloor + \lfloor 2^\lambda x_2 \rfloor \; mod \; 2$. It is a chessboard pattern based on the standard lattice $\mathbb{Z}^2 \subset \mathbb{R}^2$. Let us construct two bounded weak solutions $\zeta_\lambda$ and $\tilde{\zeta}_\lambda$ of (BVP) along $\boldsymbol{b}_\lambda$ with initial datum $\bar{\zeta}_\lambda$. As $\boldsymbol{b}_\lambda \in L^1_{loc}([0,1); BV_{loc}(\mathbb{R}^2;\mathbb{R}^2))$, there exists a unique solution of (IVP) along $\boldsymbol{b}_\lambda$ on $[0,1]$. We can now describe it. For $0 < t < 1/2$, we have $\zeta_\lambda(t,x)\mathscr{L}^d = \boldsymbol{X}_\lambda(t,0,\cdot)_\# \bar{\zeta}_\lambda(x)\mathscr{L}^d$. Using property (R), we have

$$\zeta_\lambda(1/2,x) = 1 - \bar{\zeta}_\lambda(2x). \tag{3.2}$$

By self-similarity of $\boldsymbol{b}_\lambda$, we then have $\zeta_\lambda(t,x) = \zeta_\lambda(t/2^k + 1 - 1/2^k, 2^k x)$ for $1 - 1/2^k \leq t < 1 - 1/2^{k+1}$. Using recursively the appropriately scaled version of (3.2) we can check that

$$\zeta_\lambda(1-1/2^k,x) = \bar{\zeta}_\lambda(2^k x) \quad \text{for } k \text{ even}, \qquad \zeta_\lambda(1-1/2^k,x) = 1 - \bar{\zeta}_\lambda(2^k x) \quad \text{for } k \text{ odd}.$$

Therefore $\zeta_\lambda(t,\cdot) \rightharpoonup 1/2$ as $t \uparrow 1$. Set $\zeta_\lambda(t,x) = \zeta_\lambda(2-t,x)$ for $1 < t < 2$. We also set $\tilde{\zeta}_\lambda(t,x) = \zeta_\lambda(t,x)$ for $0 < t < 1$ and $\tilde{\zeta}_\lambda(t,x) = 1/2$ for $1 \leq t < 2$. It can be directly checked that $\zeta_\lambda$ and $\tilde{\zeta}_\lambda$ are both bounded weak solution of (IVP) along $\boldsymbol{b}_\lambda$ with initial datum $\bar{\zeta}_\lambda$.

*3.2. Unperturbed truncations*

We now truncate the vector field $\boldsymbol{b}_\lambda$ around time 1 in a symmetric manner, and in an asymmetric manner. To each of these truncations will correspond a different weak solution of (IVP), which will converge, in the limit where the truncation is removed, to two distinct weak solutions of (IVP). We begin by defining the time symmetric truncation

$$\boldsymbol{b}^q_\lambda(t,x) := \begin{cases} \boldsymbol{b}_\lambda(t,x) & \text{if} \quad t \notin (1-2^{-q}, 1+2^{-q}), \\ 0 & \text{if} \quad t \in (1-2^{-q}, 1+2^{-q}). \end{cases}$$

Now, we define the time asymmetric truncation

$$\tilde{\boldsymbol{b}}^q_\lambda(t,x) := \begin{cases} \boldsymbol{b}_\lambda(t,x) & \text{if} \quad t \notin (1-2^{-q-2}, 1+2^{-q}), \\ 0 & \text{if} \quad t \in (1-2^{-q-2}, 1+2^{-q}). \end{cases}$$

Both of these vector fields lie in $L^\infty((0,2); BV_{loc} \cap L^\infty(\mathbb{R}^2;\mathbb{R}^2))$. Let us denote by $\zeta^q_\lambda$ and by $\tilde{\zeta}^q_\lambda$ the unique bounded weak solution of (IVP) along $\boldsymbol{b}^q_\lambda$ and $\tilde{\boldsymbol{b}}^q_\lambda$ respectively with initial datum $\bar{\zeta}_\lambda$. Observe that the unique bounded weak solution of (IVP) along $\boldsymbol{b}^q_\lambda$ with initial datum $\zeta_\lambda$ is

$$\zeta^q_\lambda(t,x) := \begin{cases} \zeta_\lambda(t,x) & \text{if} \quad t \notin (1-2^{-q}, 1+2^{-q}), \\ \zeta_\lambda(1-2^{-q},x) & \text{if} \quad t \in (1-2^{-q}, 1+2^{-q}). \end{cases}$$

**Lemma 3.1.** *For every $t \in [0,2]$, we have that $\zeta^q_\lambda(t,\cdot) \rightharpoonup \zeta_\lambda(t,\cdot)$ in $w^* - L^\infty_x$ and $\tilde{\zeta}^q_\lambda(t,\cdot) \rightharpoonup \tilde{\zeta}_\lambda(t,\cdot)$ in $w^* - L^\infty_x$ as $q \to +\infty$.*

**Proof.** As $\zeta_\lambda$ is in $C([0,2]; w^* - L^\infty(\mathbb{R}^2))$, we clearly have that for every $t \in [0,2]$ that $\zeta^q_\lambda(t,\cdot) \rightharpoonup \zeta_\lambda(t,\cdot)$ in $w^* - L^\infty_x$ as $q \to +\infty$. The unique bounded weak solution of (IVP) along $\tilde{\boldsymbol{b}}^q_\lambda$ with initial datum $\bar{\zeta}_\lambda$ is $\tilde{\zeta}^q_\lambda$, and given by $\tilde{\zeta}^q_\lambda(t,x) = \zeta_\lambda(t,x)$ for $0 \leq t < 1-2^{-q}$, $\tilde{\zeta}^q_\lambda(t,x) = \zeta_\lambda(1+2^{-q},x)$ for $t \in (1-2^{-q}, 1+2^{-q})$. We then clearly have that $\tilde{\zeta}^q_\lambda(t,\cdot) \rightharpoonup \tilde{\zeta}_\lambda(t,\cdot)$ in $w^* - L^\infty_x$ for every $t \in [0,1]$. Let us show that for every $t \in [1,2]$, we have $\tilde{\zeta}_\lambda(t,\cdot) \rightharpoonup \tilde{\zeta}_\lambda(t,\cdot)$ $w^* - L^\infty_x$ as $q \to +\infty$. To do so, we will first prove that $\tilde{\zeta}^q_\lambda(2,\cdot) \rightharpoonup 1/2$ in $w^* - L^\infty_x$ as $q \to +\infty$. By (O) of Section 3.1, we have for every $q' \leq q$

$$\fint_S \tilde{\zeta}^{q'}_\lambda(1+2^{-q'},x)dx = 1/2 \qquad \forall S \in S_{\lambda q}. \tag{3.3}$$

First, observe that the set

$$\mathcal{D}_\lambda := \Big\{ \sum_{j=1}^N \alpha_j \mathbb{1}_S \ : \ N \in \mathbb{N},\ \alpha_j \in \mathbb{R},\ q \in \mathbb{N},\ S \in S_{\lambda q} \Big\}$$

is dense in $L^1(\mathbb{R}^2)$. Let $\varepsilon > 0$ and let $\phi \in L^1(\mathbb{R}^2)$. Then, there exists $\psi \in \mathcal{D}_\lambda$ such that $\|\psi - \phi\|_{L^1_x} < \varepsilon$ and there exists $q_0 \in \mathbb{N}$ such that for every $q \geq q_0$, we have

$$\int_{\mathbb{R}^2} \tilde{\zeta}^q_\lambda(2,x)\psi(x)dx = \int_{\mathbb{R}^2} \frac{\psi(x)}{2}dx,$$

thanks to (3.3) and by linearity of the integral. Then, thanks to the above we have

$$\begin{aligned}
&\Big|\int_{\mathbb{R}^2} \Big(\zeta^q_\lambda(2,x)\phi(x) - \frac{\phi(x)}{2}\Big)dx\Big| \\
&= \Big|\int_{\mathbb{R}^2} \Big(\zeta^q_\lambda(2,x)\phi(x) - \zeta^q_\lambda(2,x)\psi(x) + \zeta^q_\lambda(2,x)\psi(x) - \frac{\psi(x)}{2} + \frac{\psi(x)}{2} - \frac{\phi(x)}{2}\Big)dx\Big| \quad (3.4)\\
&\leq \|\zeta^q_\lambda\|_{L^\infty_x}\|\phi - \psi\|_{L^1_x} + \frac{\|\psi - \phi\|_{L^1_x}}{2} < 3\varepsilon/2.
\end{aligned}$$

$\varepsilon$ was arbitrary, we have $\tilde{\zeta}^q_\lambda(2,\cdot) \rightharpoonup 1/2$ in $w^* - L^\infty_x$.

To conclude, we have to argue that for every $t \in [1,2]$ we have $\tilde{\zeta}^q_\lambda(t,\cdot) \rightharpoonup \tilde{\zeta}(t,\cdot)$ in $w^* - L^\infty_x$ as $q \to +\infty$. The sequence $(\tilde{\zeta}^q_\lambda)_{q\in\mathbb{N}}$ has a limit point in $C([0,2]; w^* - L^\infty(\mathbb{R}^2))$ by a compactness argument (see for instance Lemma 2.4 in [20]), and it is a solution of (BVP) along $\boldsymbol{b}_\lambda$ with boundary datum $\bar{\rho} = 1/2$ at time 2. As the constant solution is the unique bounded weak solution of (BVP) on $[1,2] \times \mathbb{R}^2$ along $\boldsymbol{b}_\lambda$ with boundary datum $1/2$ at time 2, we have that for every $t \in [1,2]$, it holds that $\tilde{\zeta}^q_\lambda(t,\cdot) \rightharpoonup 1/2$ in $w^* - L^\infty_x$. This proves the thesis. □

## 4. Perturbed vector fields

We perturb the vector fields $(\boldsymbol{b}_\lambda)_{\lambda\in\mathbb{N}}$ by an arbitrary smooth vector field in order to obtain a dense family in $L^p_t W^{s,p}_x$. Throughout this section, $\boldsymbol{w} : [0,2] \times \mathbb{R}^2 \to \mathbb{R}^2$ will always be a vector field in $L^\infty((0,2); C^2_c(\mathbb{R}^2;\mathbb{R}^2))$. We shall now denote by $\boldsymbol{X}_{\boldsymbol{w}}$ the flow starting at time 1 along $\boldsymbol{w}$.

### *4.1. The building blocks*

We then define the vector fields on $[0,2] \times \mathbb{R}^2$

$$\boldsymbol{b}_{\lambda,\boldsymbol{w}}(t,x) := D_x\boldsymbol{X}_{\boldsymbol{w}}(t, \boldsymbol{X}^{-1}_{\boldsymbol{w}}(t,x)).\boldsymbol{b}_\lambda(t, \boldsymbol{X}^{-1}_{\boldsymbol{w}}(t,x)) + \boldsymbol{w}(t,x). \quad (4.1)$$

The particular structure for this vector field is chosen so that, the flow along $\boldsymbol{b}_{\lambda,\boldsymbol{w}}$ started at time 1 will be given by the composition of the flow of the large scales $\boldsymbol{X}_{\boldsymbol{w}}$ with the flow of the small scales, namely the regular Lagrangian flow of $\boldsymbol{b}_\lambda$ started from time 1. This will explicitly state in Corollary 4.6. We then define the set

$$D := \Big\{ \boldsymbol{b}_{\lambda,\boldsymbol{w}} \ : \ \lambda \in \mathbb{N},\ \boldsymbol{w} \in L^\infty\big((0,2); C^2_c(\mathbb{R}^2;\mathbb{R}^2)\big) \Big\}. \tag{4.2}$$

**Remark 4.1.** The set $D$ is not bounded at any time in $[0,2]$ in the sense of Definition 1.2. Indeed, for every $\gamma > 0$, we can choose $\boldsymbol{w}$ in $C([0,2]; C^1_c(\mathbb{R}^2;\mathbb{R}^2))$ such that $\|\boldsymbol{w}(t,\cdot)\|_{L^\infty_x} > 2\gamma$ for every $t \in [0,2]$. Then, we can choose $\lambda \in \mathbb{N}$ large enough that for every interval $I$ with non-empty interior, we have $\|D_x \boldsymbol{X}_{\boldsymbol{w}}.\boldsymbol{b}_\lambda\|_{L^\infty(I;L^\infty_x)} \leq \gamma$, whence $\|\boldsymbol{b}_{\lambda,\boldsymbol{w}}\|_{L^\infty(I;L^\infty_x)} \geq \gamma$. This proves that $D$ is not bounded at any time in $[0,2]$.

Recall that the total variation of a vector field $\boldsymbol{u} : \mathbb{R}^d \to \mathbb{R}^d$ over an open set $\Omega \subset \mathbb{R}^d$ is given by

$$V(\boldsymbol{u},\Omega) := \sup \Big\{ \Big| \int_\Omega \boldsymbol{u}(x)\, \mathrm{div}_x\, \boldsymbol{\phi}(x) dx \Big| \ : \ \boldsymbol{\phi} \in C^1_c(\Omega,\mathbb{R}^d),\ \|\boldsymbol{\phi}\|_{L^\infty(\Omega)} \leq 1 \Big\}. \tag{4.3}$$

For any $\boldsymbol{u} \in BV_{loc}(\mathbb{R}^d;\mathbb{R}^d)$, we have that $V(\boldsymbol{u},\Omega) = |D\boldsymbol{u}|(\Omega)$, and the class of vector fields with locally bounded variation is characterised as

$$BV_{loc}(\mathbb{R}^d;\mathbb{R}^d) = \Big\{ \boldsymbol{u} \in L^1_{loc}(\mathbb{R}^d;\mathbb{R}^d) \ : \ V(\boldsymbol{u},\Omega) < +\infty \text{ for every precompact open subset } \Omega \Big\}. \tag{4.4}$$

This is for instance proven in [3, Proposition 3.6]. Given real numbers $\alpha_1, \ldots, \alpha_n$, the notation $a \lesssim_{\alpha_1,\ldots,\alpha_n} b$ will mean that there is a constant $C = C(\alpha_1,\ldots,\alpha_n)$ depending only on $\alpha_1,\ldots,\alpha_n$ such that $a \leq Cb$ for every $a$ and $b$ varying some specified sets.[1] Let us now prove that $D$ is dense.

**Lemma 4.2.** *Let $1 \leq p < +\infty$ and $0 \leq s < 1$ such that $p < 1/s$. Then the set $D$ is dense in $L^p((0,2); W^{s,p}_{loc}(\mathbb{R}^2;\mathbb{R}^2))$.*

**Proof.** Let $R > 0$ and let $\boldsymbol{w} \in C^\infty_c\big((0,2)\times\mathbb{R}^2;\mathbb{R}^2\big)$. For every $k \in \mathbb{N}$, let $I_k := (1 - 2^{-k}, 1 - 2^{-k-1}) \cup (1 + 2^{-k-1}, 1 + 2^{-k})$. Observe that the set $C^\infty_c\big((0,2)\times\mathbb{R}^2;\mathbb{R}^2\big)$ is dense in $L^p((0,2); W^{s,p}_{loc}(\mathbb{R}^2;\mathbb{R}^2))$. So it suffices to prove that the sequence $(\boldsymbol{b}_{\lambda,\boldsymbol{w}})_{\lambda\in\mathbb{N}}$ converges to $\boldsymbol{w}$ in $L^p((0,2); W^{s,p}_{loc}(\mathbb{R}^2;\mathbb{R}^2))$. Let $0 < \sigma < 1$ be a parameter to be determined later.

**Step 1** ($W^{\sigma,1}_{loc}$ control.) By Gagliardo-Nirenberg interpolation (see [8]), for every $t \in [0,2]$ we have

$$\begin{aligned} &\|\boldsymbol{b}_{\lambda,\boldsymbol{w}}(t,\cdot) - \boldsymbol{w}(t,\cdot)\|_{W^{\sigma,1}(B_R(0)))} \\ &\quad \lesssim_{R,\sigma} \|\boldsymbol{b}_{\lambda,\boldsymbol{w}}(t,\cdot) - \boldsymbol{w}(t,\cdot)\|^{(1-\sigma)}_{L^1(B_R(0))} \|\boldsymbol{b}_{\lambda,\boldsymbol{w}}(t,\cdot) - \boldsymbol{w}(t,\cdot)\|^\sigma_{BV(B_R(0))}. \end{aligned} \tag{4.5}$$

[1] The notation $a \lesssim b$ will mean that there is a numerical constant $C > 0$ such that $a \leq Cb$ for every $a$ and $b$ varying in some specified sets.

Let $t \in I_k$, and let us estimate $\|\boldsymbol{b}_{\lambda,\boldsymbol{w}}(t,\cdot) - \boldsymbol{w}(t,\cdot)\|_{BV(B_R(0))}$. Let $N \in \mathbb{N}$ be sufficiently large that $B_{R+2\|\boldsymbol{w}\|_{L^\infty_{t,x}}}(0) \subset (-N,N)^2$. Denote the constant $E := \exp\left(\int_0^2 \|\boldsymbol{w}(u,\cdot)\|_{C^1_x} du\right)$. For every $\boldsymbol{\phi} \in C^1_c(B_R(0);\mathbb{R}^2)$ with $\|\boldsymbol{\phi}\|_{L^\infty(B_R(0))} \le 1$, we then have by a change of variable, and estimating using subsection 2.3 that

$$\begin{aligned}
&\Big|\int_{B_R(0)} D_x\boldsymbol{X}_{\boldsymbol{w}}(t,\boldsymbol{X}_{\boldsymbol{w}}^{-1}(t,x)).\boldsymbol{b}_\lambda(t,\boldsymbol{X}_{\boldsymbol{w}}^{-1}(t,x))\,\mathrm{div}_x\,\boldsymbol{\phi}(x)dx\Big| \\
&\le \Big|\int_{(-N,N)^2} D_x\boldsymbol{X}_{\boldsymbol{w}}(t,y).\boldsymbol{b}_\lambda(t,y)\,\mathrm{div}_x\,\boldsymbol{\phi}(\boldsymbol{X}_{\boldsymbol{w}}(t,y))J\boldsymbol{X}_{\boldsymbol{w}}(t,y)dy\Big| \\
&\lesssim E^2\Big|\int_{(-N,N)^2} \boldsymbol{b}_\lambda(t,y)\,\mathrm{div}_x\,\boldsymbol{\phi}(\boldsymbol{X}_{\boldsymbol{w}}(t,y))dy\Big| \\
&\lesssim E^2\|\boldsymbol{b}_\lambda(t,\cdot)\|_{BV((-N,N)^2)} \\
&\lesssim E^2\|\boldsymbol{b}_0(t,\cdot)\|_{BV((-2N,2N)^2)}.
\end{aligned} \tag{4.6}$$

In the second to last inequality, we have used that $\boldsymbol{\phi} \circ \boldsymbol{X}_{\boldsymbol{w}}(t,\cdot) \in C^1_c((-N,N)^2;\mathbb{R}^2)$ and $\|\boldsymbol{\phi} \circ \boldsymbol{X}_{\boldsymbol{w}}(t,\cdot)\|_{L^\infty((-N,N)^2)} \le 1$, and in the last inequality that $\|\boldsymbol{b}_\lambda(t,\cdot)\|_{BV((-N,N)^2)} \le \|\boldsymbol{b}_0(t,\cdot)\|_{BV((-2N,2N)^2)}$. Therefore,

$$\|\boldsymbol{b}_{\lambda,\boldsymbol{w}} - \boldsymbol{w}\|_{BV(B_R(0))} \lesssim E^2\|\boldsymbol{b}_0(t,\cdot)\|_{BV((-2N,2N)^2)}. \tag{4.7}$$

Furthermore, we have by a change of variable and then estimating using subsection 2.3 that

$$\begin{aligned}
&\int_{B_R(0)} |D_x\boldsymbol{X}_{\boldsymbol{w}}(t,\boldsymbol{X}_{\boldsymbol{w}}^{-1}(t,x)).\boldsymbol{b}_\lambda(t,\boldsymbol{X}_{\boldsymbol{w}}^{-1}(t,x))|dx \\
&\lesssim \int_{(-N,N)^2} |D_x\boldsymbol{X}_{\boldsymbol{w}}(t,y).\boldsymbol{b}_\lambda(t,y)J\boldsymbol{X}_{\boldsymbol{w}}(t,y)|dy \\
&\lesssim E^2\|\boldsymbol{b}_\lambda\|_{L^1((-N,N)^2)} \\
&\lesssim 2^{-\lambda}E^2\|\boldsymbol{b}_0\|_{L^1((-2N,2N)^2)}.
\end{aligned} \tag{4.8}$$

Thus,

$$\|\boldsymbol{b}_{\lambda,\boldsymbol{w}} - \boldsymbol{w}\|_{L^1(B_R(0))} \lesssim 2^{-\lambda}E^2\|\boldsymbol{b}_0(t,\cdot)\|_{L^1((-2N,2N))^2)}. \tag{4.9}$$

Now, observe that there exists a constant $C_N > 0$, which depends on $N$ only, and such that for every $k \in \mathbb{N}$, and every $t \in I_k$, we have

$$\|\boldsymbol{b}_0(t,\cdot)\|_{L^1((-2N,2N)^2)} \le C_N \qquad \text{and} \qquad \|\boldsymbol{b}_0(t,\cdot)\|_{BV((-2N,2N)^2)} \le C_N 2^k. \tag{4.10}$$

Gathering (4.5), (4.7), (4.9) and (4.10), we have

$$\begin{aligned}\|\boldsymbol{b}_{\lambda,\boldsymbol{w}}(t,\cdot)-\boldsymbol{w}(t,\cdot)\|_{W^{\sigma,1}(B_R(0))} &\lesssim_{N,R,\sigma} 2^{-\sigma\lambda}E^2\|\boldsymbol{b}_0(t,\cdot)\|^{\sigma}_{BV((-2N,2N)^2)} \\ &\lesssim_{N,R,\sigma} 2^{-\sigma\lambda}E^2 2^{k\sigma}.\end{aligned} \tag{4.11}$$

**Step 2** ($L^p((0,2);W^{s,p}_{loc})$ control.) Now fix the parameter $\sigma = ps$. By Gagliardo-Nirenberg interpolation, for every $t\in[0,2]$ we have

$$\begin{aligned}&\|\boldsymbol{b}_{\lambda,\boldsymbol{w}}(t,\cdot)-\boldsymbol{w}(t,\cdot)\|_{W^{s,p}(B_R(0))} \\ &\quad\lesssim_{R,p,s} \|\boldsymbol{b}_{\lambda,\boldsymbol{w}}(t,\cdot)-\boldsymbol{w}(t,\cdot)\|^{1-s/\sigma}_{L^\infty(B_R(0))}\|\boldsymbol{b}_{\lambda,\boldsymbol{w}}(t,\cdot)-\boldsymbol{w}(t,\cdot)\|^{s/\sigma}_{W^{\sigma,1}(B_R(0))}.\end{aligned} \tag{4.12}$$

Note also that

$$\|\boldsymbol{b}_{\lambda,\boldsymbol{w}}(t,\cdot)-\boldsymbol{w}(t,\cdot)\|_{L^\infty(\mathbb{R}^d)} \le 2^{-\lambda}E\|\boldsymbol{b}_0\|_{L^\infty((0,2)\times\mathbb{R}^2)}.$$

Therefore using (4.11), for every $t\in I_k$, we have

$$\|\boldsymbol{b}_{\lambda,\boldsymbol{w}}(t,\cdot)-\boldsymbol{w}(t,\cdot)\|_{W^{s,p}(B_R(0))} \lesssim_{N,R,p,s} 2^{-\lambda(1-s(1-1/\sigma))}E^2 2^{ks}. \tag{4.13}$$

So, we have

$$\|\boldsymbol{b}_{\lambda,\boldsymbol{w}}(t,\cdot)-\boldsymbol{w}(t,\cdot)\|_{L^p((0,2);W^{s,p}(B_R(0)))} \lesssim_{N,R,p,s} 2^{-\lambda(1-s(1-1/\sigma))}E^2\sum_{k\in\mathbb{N}}2^{ks(s-1/p)}. \tag{4.14}$$

The geometric series is finite because $s-1/p<0$. As $1-s(1-1/\sigma)>0$ and as $R$ was arbitrary, it follows that $\boldsymbol{b}_{\lambda,\boldsymbol{w}}$ converges to $\boldsymbol{w}$ in $L^p((0,2);W^{s,p}_{loc})$ as $\lambda\to+\infty$. This proves the thesis. □

### *4.2. Perturbed truncations*

We now truncate the rough part of vector fields $\boldsymbol{b}_{\lambda,\boldsymbol{w}}$ at time 1 in a symmetric manner, and in an asymmetric manner. We will conclude this section by exhibiting a composition rule for their regular Lagrangian flow starting at time 1, thereby decoupling the dynamics from the smooth and the non-smooth part of the vector field.

We define the symmetric truncation

$$\boldsymbol{b}^q_{\lambda,\boldsymbol{w}}(t,x) := D_x\boldsymbol{X}_{\boldsymbol{w}}(t,\boldsymbol{X}^{-1}_{\boldsymbol{w}}(t,x)).\boldsymbol{b}^q_{\lambda}(t,\boldsymbol{X}^{-1}_{\boldsymbol{w}}(t,x))+\boldsymbol{w}(t,x), \tag{4.15}$$

and the asymmetric truncation

$$\tilde{\boldsymbol{b}}^q_{\lambda,\boldsymbol{w}}(t,x) := D_x\boldsymbol{X}_{\boldsymbol{w}}(t,\boldsymbol{X}^{-1}_{\boldsymbol{w}}(t,x)).\tilde{\boldsymbol{b}}^q_{\lambda}(t,\boldsymbol{X}^{-1}_{\boldsymbol{w}}(t,x))+\boldsymbol{w}(t,x). \tag{4.16}$$

Notice that the divergence of both of these vector fields may be singular with respect to $\mathscr{L}^2$ in contrast with the work of De Lellis and Giri [12].

**Lemma 4.3.** *The vector fields* $\boldsymbol{b}^q_{\lambda,\boldsymbol{w}}$ *and* $\tilde{\boldsymbol{b}}^q_{\lambda,\boldsymbol{w}}$ *belong to* $L^\infty((0,2); BV_{loc} \cap L^\infty(\mathbb{R}^2;\mathbb{R}^2))$.

**Proof.** Clearly, it suffices to prove the thesis for $\tilde{\boldsymbol{b}}^q_{\lambda,\boldsymbol{w}}$. It is readily seen that $\tilde{\boldsymbol{b}}^q_{\lambda,\boldsymbol{w}}$ is in $L^\infty((0,2)\times\mathbb{R}^2;\mathbb{R}^2)$. To prove the thesis, it is enough to show that

$$L^1((0,2)) \ni \psi \longmapsto \int_0^2 V(\tilde{\boldsymbol{b}}^q_{\lambda,\boldsymbol{w}}(t,\cdot),\Omega)\psi(t)dt,$$

is a bounded linear functional for every precompact $\Omega \subset \mathbb{R}^2$. Let $\Omega$ be a precompact subset of $\mathbb{R}^2$, let $\psi \in L^1((0,1))$, and let $\boldsymbol{\phi} \in C^1_c(\Omega;\mathbb{R}^2)$. Let $R > 0$ be sufficiently large that $\boldsymbol{X}^{-1}_{\boldsymbol{w}}(t,\Omega) \subset B_R(0)$ for every $t \in [0,2]$. Notice that such an $R$ exists by finite speed of propagation for a bounded smooth vector field. We then have

$$\begin{aligned}
&\int_0^2 \Big| \int_{\mathbb{R}^2} \tilde{\boldsymbol{b}}^q_{\lambda,\boldsymbol{w}}(t,x)\,\mathrm{div}_x\,\boldsymbol{\phi}(x)dx \Big| \psi(t)dt \\
&\le \underbrace{\int_0^2 \Big| \int_{\mathbb{R}^2} D_x\boldsymbol{X}_{\boldsymbol{w}}(t,\boldsymbol{X}^{-1}_{\boldsymbol{w}}(t,x)).\boldsymbol{b}^q_\lambda(t,\boldsymbol{X}^{-1}_{\boldsymbol{w}}(t,x))\,\mathrm{div}_x\,\boldsymbol{\phi}(x)dx \Big| \psi(t)dt}_{I} \\
&+ \underbrace{\int_0^2 \Big| \int_{\mathbb{R}^2} \boldsymbol{w}(t,x)\,\mathrm{div}_x\,\boldsymbol{\phi}(x)dx \Big| \psi(t)dt}_{II}.
\end{aligned} \tag{4.17}$$

For the second term, performing an integration by parts gives that

$$II = \int_0^2 \Big| \int_{\mathbb{R}^2} D_x\boldsymbol{w}(t,x).\boldsymbol{\phi}(x)dx \Big| \psi(t)dt \le \|\boldsymbol{w}\|_{C^1_{t,x}} \|\boldsymbol{\phi}\|_{L^\infty_x} \|\psi\|_{L^1_t}.$$

So it suffices to treat the first term. A change of variables $y = \boldsymbol{X}^{-1}_{\boldsymbol{w}}(t,x)$, and the estimates of subsection 2.3 then give

$$I = \int_0^2 \Big| \int_{B_R(0)} D_x \boldsymbol{X}_{\boldsymbol{w}}(t,y).\boldsymbol{b}_\lambda^q(t,y) J\boldsymbol{X}_{\boldsymbol{w}}(t,y) \operatorname{div}_x \boldsymbol{\phi}(\boldsymbol{X}_{\boldsymbol{w}}(t,y))dy \Big| \psi(t)dt$$

$$\leq \exp\Big(2\int_0^2 \|\boldsymbol{w}(u,\cdot)\|_{C_x^1}\Big) \int_0^2 V(\boldsymbol{b}_\lambda^q(t,\cdot), B_R(0))\psi(t)dt$$

$$\leq \exp\Big(2\int_0^2 \|\boldsymbol{w}(u,\cdot)\|_{C_x^1}\Big) \|\boldsymbol{b}_\lambda^q\|_{L^\infty((0,2);BV(B_R(0)))} \|\psi\|_{L^1((0,2))}.$$

As $\boldsymbol{\phi}$ was arbitrary in $C_c^1(\Omega)$, there exists a constant $C > 0$ independant of $\psi$ such that

$$\int_0^2 V(\tilde{\boldsymbol{b}}_{\lambda,\boldsymbol{w}}^q(t,\cdot),\Omega)\psi(t)dt \leq C\|\psi\|_{L_t^1}. \tag{4.18}$$

As $\psi$ was arbitrary in $L^1((0,2))$ and $\Omega$ was arbitrary, the thesis follows. $\square$

Although $\operatorname{div}_x \boldsymbol{b}_\lambda(t,\cdot) = 0$ for a.e. $t \in [0,2]$, and $\operatorname{div}_x \boldsymbol{w}(t,\cdot)$ admits a bounded density with respect to $\mathscr{L}^2$ for every $t \in [0,2]$, the nonlinear transform performed on $\boldsymbol{b}_\lambda^q$ and $\tilde{\boldsymbol{b}}_\lambda^q$ implies that $\operatorname{div}_x \boldsymbol{b}_{\lambda,\boldsymbol{w}}^q(t,\cdot)$ and $\operatorname{div}_x \tilde{\boldsymbol{b}}_{\lambda,\boldsymbol{w}}^q(t,\cdot)$ may have non-vanishing singular part with respect to $\mathscr{L}^2$. Thanks to Theorem 2.5, it is nonetheless possible to identify the regular Lagrangian flows of both $\boldsymbol{b}_{\lambda,\boldsymbol{w}}^q$ and $\tilde{\boldsymbol{b}}_{\lambda,\boldsymbol{w}}^q$ starting at 1 as both these vector fields will be shown to be nearly incompressible. We shall first begin by proving a chain rule. We denote by $\boldsymbol{X}_{\boldsymbol{b}_\lambda^q}$ the regular Lagrangian flows starting at time 1 along $\boldsymbol{b}_\lambda^q$ and by $\boldsymbol{X}_{\tilde{\boldsymbol{b}}_\lambda^q}$ the regular Lagrangian flows starting at time 1 along $\tilde{\boldsymbol{b}}_\lambda^q$. We then have the following chain rule.

**Lemma 4.4.** *For every $q \in \mathbb{N}$, the following chain rule identities hold for $\mathscr{L}^3$-a.e. $(t,x) \in (0,2) \times \mathbb{R}^2$*

$$\begin{aligned}
(i)\ \ \frac{d}{dt}\boldsymbol{X}_{\boldsymbol{w}}(t,\boldsymbol{X}_{\boldsymbol{b}_\lambda^q}(t,x)) =& D_x\boldsymbol{X}_{\boldsymbol{w}}(t,\boldsymbol{X}_{\boldsymbol{b}_\lambda^q}(t,x)).\boldsymbol{b}_\lambda^q(t,\boldsymbol{X}_{\boldsymbol{b}_\lambda^q}(t,x)) \\
&+ \boldsymbol{w}(t,\boldsymbol{X}_{\boldsymbol{w}}(t,\boldsymbol{X}_{\boldsymbol{b}_\lambda^q}(t,x))), \\
(ii)\ \ \frac{d}{dt}\boldsymbol{X}_{\boldsymbol{w}}(t,\boldsymbol{X}_{\tilde{\boldsymbol{b}}_\lambda^q}(t,x)) =& D_x\boldsymbol{X}_{\boldsymbol{w}}(t,\boldsymbol{X}_{\tilde{\boldsymbol{b}}_\lambda^q}(t,x)).\tilde{\boldsymbol{b}}_\lambda^q(t,\boldsymbol{X}_{\tilde{\boldsymbol{b}}_\lambda^q}(t,x)) \\
&+ \boldsymbol{w}(t,\boldsymbol{X}_{\boldsymbol{w}}(t,\boldsymbol{X}_{\tilde{\boldsymbol{b}}_\lambda^q}(t,x))).
\end{aligned} \tag{4.19}$$

**Proof.** Let us prove $(i)$. The proof of $(ii)$ is completely analogous. Let $\eta \in C_c^\infty(\mathbb{R})$ be a standard mollifier, and denote $\eta^k(t) = k\eta(kt)$. Let $\phi \in C_c^\infty((0,2))$ be an arbitrary test function. So that time convolution will be well-defined, we extend $\boldsymbol{X}_{\boldsymbol{w}}(t,\cdot)$ and $\boldsymbol{X}_{\boldsymbol{b}_\lambda^q}(t,\cdot)$ to $t \in \mathbb{R}$ by setting $\boldsymbol{X}_{\boldsymbol{w}}(t,x) = \boldsymbol{X}_{\boldsymbol{w}}(0,x)$ and $\boldsymbol{X}_{\boldsymbol{b}_\lambda^q}(t,x) = \boldsymbol{X}_{\boldsymbol{b}_\lambda^q}(0,x)$ for $t \leq 0$, and setting $\boldsymbol{X}_{\boldsymbol{w}}(t,x) = \boldsymbol{X}_{\boldsymbol{w}}(1,x)$ and $\boldsymbol{X}_{\boldsymbol{b}^q_\lambda}(t,x) = \boldsymbol{X}_{\boldsymbol{b}^q_\lambda}(1,x)$ for $t \geq 2$ and a.e. $x \in \mathbb{R}^2$. We then have for a.e. $x \in \mathbb{R}^2$

$$
\begin{aligned}
&\int_0^2 \boldsymbol{X}_{\boldsymbol{w}}(t, \boldsymbol{X}_{\boldsymbol{b}^q_\lambda}(t,x))\partial_t\phi(t)dt \\
&\overset{1}{=} \lim_{k\to+\infty} \int_0^2 \eta^k \star_t \boldsymbol{X}_{\boldsymbol{w}}(t, \eta^k \star_t \boldsymbol{X}_{\boldsymbol{b}^q_\lambda}(t,x))\partial_t\phi(t)dt \\
&\overset{2}{=} -\lim_{k\to+\infty} \int_0^2 \frac{d}{dt}\Big[\eta^k \star_t \boldsymbol{X}_{\boldsymbol{w}}(t, \eta^k \star_t \boldsymbol{X}_{\boldsymbol{b}^q_\lambda}(t,x))\Big]\phi(t)dt \\
&\overset{3}{=} -\lim_{k\to+\infty} \int_0^2 \Big[\eta^k \star_t \partial_t\boldsymbol{X}_{\boldsymbol{w}}(t, \eta^k \star_t \boldsymbol{X}_{\boldsymbol{b}^q_\lambda}(t,x)) \\
&\quad + D_x\boldsymbol{X}_{\boldsymbol{w}}(t, \eta^k \star_t \boldsymbol{X}_{\boldsymbol{b}^q_\lambda}(t,x)).\eta^k \star_t \partial_t\boldsymbol{X}_{\boldsymbol{b}^q_\lambda}(t,x)\Big]\phi(t)dt \\
&\overset{4}{=} -\lim_{k\to+\infty} \int_0^2 \Big[\eta^k \star_t \boldsymbol{w}(t, \boldsymbol{X}_{\boldsymbol{w}}(t, \eta^k \star_t \boldsymbol{X}_{\boldsymbol{b}^q_\lambda}(t,x)) \\
&\quad + D_x\boldsymbol{X}_{\boldsymbol{w}}(t, \eta^k \star_t \boldsymbol{X}_{\boldsymbol{b}^q_\lambda}(t,x)).\eta^k \star_t \boldsymbol{b}^q_\lambda(t,x)\Big]\phi(t)dt \\
&\overset{5}{=} -\int_0^2 \Big[\boldsymbol{w}(t, \boldsymbol{X}_{\boldsymbol{w}}(t, \boldsymbol{X}_{\boldsymbol{b}^q_\lambda}(t,x)) + D_x\boldsymbol{X}_{\boldsymbol{w}}(t, \boldsymbol{X}_{\boldsymbol{b}^q_\lambda}(t,x)).\boldsymbol{b}^q_\lambda(t,x)\Big]\phi(t)dt,
\end{aligned}
\tag{4.20}
$$

where equality 1 follows by dominated convergence, equality 2 by an intergration by parts, equality 3 by the classical chain rule, equality 4 holds for a.e. $x \in \mathbb{R}^2$ since $\boldsymbol{X}_{\boldsymbol{w}}$ is a flow along $\boldsymbol{w}$ and $\boldsymbol{X}_{\boldsymbol{b}^q_\lambda}$ is a regular Lagrangian flow along $\boldsymbol{b}^q_\lambda$, equality 5 holds by dominated convergence. As $\phi$ was arbitrary, the thesis follows. □

Using the above chain rule, we can prove the following lemma.

**Lemma 4.5.** *The vector fields* $\boldsymbol{b}^q_{\lambda,\boldsymbol{w}}$ *and* $\tilde{\boldsymbol{b}}^q_{\lambda,\boldsymbol{w}}$ *are nearly incompressible.*

**Proof.** Let us begin by showing that there exists a regular Lagrangian flow $\boldsymbol{X}_{\boldsymbol{b}^q_{\lambda,\boldsymbol{w}}}$ starting at 1 along $\boldsymbol{b}^q_{\lambda,\boldsymbol{w}}$ for some compressibility constant $C > 0$. We set $\boldsymbol{X}_{\boldsymbol{b}^q_{\lambda,\boldsymbol{w}}}(t,x) := \boldsymbol{X}_{\boldsymbol{w}}(t, \boldsymbol{X}_{\boldsymbol{b}^q_\lambda}(t,x))$ for $t \in [0,2]$. Using Lemma 4.4, we have for $\mathscr{L}^3$-a.e. $(t,x) \in (0,2) \times \mathbb{R}^2$ that

$$\begin{aligned}\frac{d}{dt}\boldsymbol{X}_{\boldsymbol{b}^q_{\lambda,\boldsymbol{w}}}(t,x) &= \frac{d}{dt}\big(\boldsymbol{X}_{\boldsymbol{w}}(t,\boldsymbol{X}_{\boldsymbol{b}^q_\lambda}(t,x))\big)\\ &= D_x\boldsymbol{X}_{\boldsymbol{w}}(t,\boldsymbol{X}_{\boldsymbol{b}^q_\lambda}(t,x)).\boldsymbol{b}^q_\lambda(t,\boldsymbol{X}_{\boldsymbol{b}^q_\lambda}(t,x))\\ &\qquad + \boldsymbol{w}(t,\boldsymbol{X}_{\boldsymbol{w}}(t,\boldsymbol{X}_{\boldsymbol{b}^q_\lambda}(t,x)))\\ &= \boldsymbol{b}^q_{\lambda,\boldsymbol{w}}(t,\boldsymbol{X}_{\boldsymbol{w}}(t,\boldsymbol{X}_{\boldsymbol{b}^q_\lambda}(t,x)))\\ &= \boldsymbol{b}^q_{\lambda,\boldsymbol{w}}(t,\boldsymbol{X}_{\boldsymbol{b}^q_{\lambda,\boldsymbol{w}}}(t,x)),\end{aligned} \tag{4.21}$$

where in the second to last equality, we have used (4.15) and in the last equality, we have used the definition of $\boldsymbol{X}_{\boldsymbol{b}^q_{\lambda,\boldsymbol{w}}}$. Therefore, for a.e. $x \in \mathbb{R}^2$, we have that $[0,2] \ni t \mapsto \boldsymbol{X}_{\boldsymbol{b}^q_{\lambda,\boldsymbol{w}}}(t,x)$ is an absolutely continuous integral curve of $\boldsymbol{b}^q_{\lambda,\boldsymbol{w}}$ satisfying $\boldsymbol{X}_{\boldsymbol{b}^q_{\lambda,\boldsymbol{w}}}(1,x) = x$.

Observe also that for $t \in [0,2]$, we have

$$\begin{aligned}\boldsymbol{X}_{\boldsymbol{b}^q_{\lambda,\boldsymbol{w}}}(t,\cdot)_{\#}\mathscr{L}^d &= \boldsymbol{X}_{\boldsymbol{w}}(t,\boldsymbol{X}_{\boldsymbol{b}^q_\lambda}(t,\cdot))_{\#}\mathscr{L}^d\\ &= \boldsymbol{X}_{\boldsymbol{w}}(t,\cdot)_{\#}\boldsymbol{X}_{\boldsymbol{b}^q_\lambda}(t,\cdot)_{\#}\mathscr{L}^d\\ &= \boldsymbol{X}_{\boldsymbol{w}}(t,\cdot)_{\#}\mathscr{L}^d.\end{aligned} \tag{4.22}$$

Therefore, we have

$$\begin{aligned}\exp\Big(-\int_0^2 \|[\operatorname{div}_x \boldsymbol{w}(s,\cdot)]^+\|_{L^\infty_x}\,ds\Big)\mathscr{L}^d &\le \boldsymbol{X}_{\boldsymbol{b}^q_{\lambda,\boldsymbol{w}}}(t,\cdot)_{\#}\mathscr{L}^d\\ &\le \exp\Big(\int_0^2 \|[\operatorname{div}_x \boldsymbol{w}(s,\cdot)]^-\|_{L^\infty_x}\,ds\Big)\mathscr{L}^d.\end{aligned} \tag{4.23}$$

Thus, we can set

$$C := \exp\Big(\int_0^2 \|\operatorname{div}_x \boldsymbol{w}(s,\cdot)\|_{L^\infty_x}\,ds\Big),$$

as a compressibility constant of $\boldsymbol{X}_{\boldsymbol{b}^q_{\lambda,\boldsymbol{w}}}$. Therefore, the density

$$\rho^q(t,\cdot)\mathscr{L}^d = \boldsymbol{X}_{\boldsymbol{b}^q_{\lambda,\boldsymbol{w}}}(t,\cdot)_{\#}\mathscr{L}^d,$$

solves (BVP) along $\boldsymbol{b}^q_{\lambda,\boldsymbol{w}}$ and satisfies $C^{-1} < \rho^q(t,x) < C$ for every $t \in [0,2]$ and for a.e. $x \in \mathbb{R}^2$. This shows that $\boldsymbol{b}^q_{\lambda,\boldsymbol{w}}$ is nearly incompressible. An entirely similar analysis also proves that $\tilde{\boldsymbol{b}}^q_{\lambda,\boldsymbol{w}}$ is nearly incompressible. This concludes the proof. □

So we have seen in this proof that the deformation of the Lebesgue measure by the flow of $\boldsymbol{b}^q_{\lambda,\boldsymbol{w}}$ only depends on $\boldsymbol{w}$. As a direct corollary of the above argument, and in view of Theorem 2.5, we have the following.

**Corollary 4.6.** *The unique regular Lagrangian flow started from* 1 *of* $\boldsymbol{b}^q_{\lambda,\boldsymbol{w}}$ *is given by* $\boldsymbol{X}_{\boldsymbol{b}^q_{\lambda,\boldsymbol{w}}}(t,x) = \boldsymbol{X}_{\boldsymbol{w}}(t, \boldsymbol{X}_{\boldsymbol{b}^q_\lambda}(t,x))$ *for a.e.* $x \in \mathbb{R}^2$, *and the unique regular Lagrangian flow started from* 1 *of* $\tilde{\boldsymbol{b}}^q_{\lambda,\boldsymbol{w}}$ *is given by* $\boldsymbol{X}_{\tilde{\boldsymbol{b}}^q_{\lambda,\boldsymbol{w}}}(t,x) = \boldsymbol{X}_{\boldsymbol{w}}(t, \boldsymbol{X}_{\tilde{\boldsymbol{b}}^q_\lambda}(t,x))$ *for a.e.* $x \in \mathbb{R}^2$.

*4.3. The weak solutions*

In the previous section, we have proven that $\boldsymbol{b}^q_{\lambda,\boldsymbol{w}}$ and $\tilde{\boldsymbol{b}}^q_{\lambda,\boldsymbol{w}}$ are nearly incompressible vector fields in $L^\infty((0,2); BV_{loc} \cap L^\infty(\mathbb{R}^2;\mathbb{R}^2))$. Therefore, by Theorem 2.3, there exist unique solutions $\rho^q_{\lambda,\boldsymbol{w}}$ and $\tilde{\rho}^q_{\lambda,\boldsymbol{w}}$ along $\boldsymbol{b}^q_{\lambda,\boldsymbol{w}}$ and $\tilde{\boldsymbol{b}}^q_{\lambda,\boldsymbol{w}}$ respectively with initial datum $\bar{\rho}$ given by

$$\bar{\rho}\mathscr{L}^2 = \boldsymbol{X}_{\boldsymbol{w}}(0,\cdot)_{\#}\bar{\zeta}_\lambda \mathscr{L}^2. \tag{4.24}$$

We now describe these unique bounded weak solutions. This is the content of the following proposition.

**Proposition 4.7.** *We have that for every* $t \in [0,2]$

$$\rho^q_{\lambda,\boldsymbol{w}}(t,\cdot)\mathscr{L}^2 = \boldsymbol{X}_{\boldsymbol{w}}(t,\cdot)_{\#}\zeta^q_\lambda(t,\cdot)\mathscr{L}^2 \quad \textit{and} \quad \tilde{\rho}^q_{\lambda,\boldsymbol{w}}(t,\cdot)\mathscr{L}^2 = \boldsymbol{X}_{\boldsymbol{w}}(t,\cdot)_{\#}\tilde{\zeta}^q_\lambda(t,\cdot)\mathscr{L}^2.$$

**Proof.** By the construction of Section 3.1, we have

$$\boldsymbol{X}_{\boldsymbol{b}^q_\lambda}(0,\cdot)_{\#}\zeta_\lambda(1-2^{-q},\cdot)\mathscr{L}^2 = \bar{\zeta}_\lambda\mathscr{L}^2 \qquad \text{and} \qquad \boldsymbol{X}_{\tilde{\boldsymbol{b}}^q_\lambda}(0,\cdot)_{\#}\zeta_\lambda(1-2^{-q-2},\cdot)\mathscr{L}^2 = \bar{\zeta}_\lambda\mathscr{L}^2.$$

Therefore, in view of Corollary 4.6, we also have

$$\begin{aligned} \boldsymbol{X}_{\boldsymbol{b}^q_{\lambda,\boldsymbol{w}}}(0,\cdot)_{\#}\zeta_\lambda(1-2^{-q},\cdot)\mathscr{L}^2 &= \boldsymbol{X}_{\boldsymbol{w}}(0,\cdot)_{\#}\boldsymbol{X}_{\boldsymbol{b}^q_\lambda}(0,\cdot)_{\#}\zeta^q_\lambda(1-2^{-q},\cdot)\mathscr{L}^2 \\ &= \boldsymbol{X}_{\boldsymbol{w}}(0,\cdot)_{\#}\bar{\zeta}_\lambda\mathscr{L}^2 = \bar{\rho}\mathscr{L}^2, \end{aligned}$$

as well as

$$\begin{aligned} \boldsymbol{X}_{\tilde{\boldsymbol{b}}^q_{\lambda,\boldsymbol{w}}}(0,\cdot)_{\#}\zeta_\lambda(1-2^{-q-2},\cdot)\mathscr{L}^2 &= \boldsymbol{X}_{\boldsymbol{w}}(0,\cdot)_{\#}\boldsymbol{X}_{\tilde{\boldsymbol{b}}^q_\lambda}(0,\cdot)_{\#}\zeta^q_\lambda(1-2^{-q-2},\cdot)\mathscr{L}^2 \\ &= \boldsymbol{X}_{\boldsymbol{w}}(0,\cdot)_{\#}\bar{\zeta}_\lambda\mathscr{L}^2 = \bar{\rho}\mathscr{L}^2. \end{aligned}$$

Therefore, the unique bounded weak solution of (BVP) along $\boldsymbol{b}^q_{\lambda,\boldsymbol{w}}$ with boundary datum $\zeta^q(1-2^{-q},\cdot)$, and the unique bounded weak solution $\rho^q_{\lambda,\boldsymbol{w}}$ of (IVP) along $\boldsymbol{b}^q_{\lambda,\boldsymbol{w}}$ with initial datum $\bar{\rho}$ coincides. Therefore, for every $t \in [0,2]$, we have

$$\begin{aligned} \rho^q_{\lambda,\boldsymbol{w}}(t,\cdot)\mathscr{L}^2 &= \boldsymbol{X}_{\boldsymbol{b}^q_{\lambda,\boldsymbol{w}}}(t,\cdot)_{\#}\zeta_\lambda(1-2^{-q},\cdot)\mathscr{L}^2 \\ &= \boldsymbol{X}_{\boldsymbol{w}}(t,\cdot)_{\#}\boldsymbol{X}_{\boldsymbol{b}^q_\lambda}(t,\cdot)_{\#}\zeta_\lambda(1-2^{-q},\cdot)\mathscr{L}^2 \\ &= \boldsymbol{X}_{\boldsymbol{w}}(t,\cdot)_{\#}\zeta^q_\lambda(t,\cdot)\mathscr{L}^2, \end{aligned}$$

where the second equality follows from Corollary 4.6. Also the unique bounded weak solution of (BVP) along $\tilde{\boldsymbol{b}}^q_{\lambda,\boldsymbol{w}}$ with boundary datum $\zeta^q(1-2^{-q-2},\cdot)$, and the unique bounded weak solution $\tilde{\rho}^q_{\lambda,\boldsymbol{w}}$ of (IVP) along $\tilde{\boldsymbol{b}}^q_{\lambda,\boldsymbol{w}}$ with initial datum $\bar{\rho}$ coincides. Therefore, for every $t \in [0,2]$, we have

$$\begin{aligned}\tilde{\rho}^q_{\lambda,\boldsymbol{w}}(t,\cdot)\mathscr{L}^2 &= \boldsymbol{X}_{\tilde{\boldsymbol{b}}^q_{\lambda,\boldsymbol{w}}}(t,\cdot)_{\#}\zeta_\lambda(1-2^{-q-2},\cdot)\mathscr{L}^2 \\ &= \boldsymbol{X}_{\boldsymbol{w}}(t,\cdot)_{\#}\boldsymbol{X}_{\tilde{\boldsymbol{b}}^q_\lambda}(t,\cdot)_{\#}\zeta_\lambda(1-2^{-q-2},\cdot)\mathscr{L}^2 \\ &= \boldsymbol{X}_{\boldsymbol{w}}(t,\cdot)_{\#}\tilde{\zeta}^q_\lambda(t,\cdot)\mathscr{L}^2,\end{aligned}$$

where the second equality follows from Corollary 4.6. The thesis follows. □

We finally set for every $t \in [0,2]$

$$\rho_{\lambda,\boldsymbol{w}}(t,\cdot)\mathscr{L}^2 = \boldsymbol{X}_{\boldsymbol{w}}(t,\cdot)_{\#}\zeta_\lambda(t,\cdot)\mathscr{L}^2, \tag{4.25}$$

as well as

$$\tilde{\rho}_{\lambda,\boldsymbol{w}}(t,\cdot)\mathscr{L}^2 = \boldsymbol{X}_{\boldsymbol{w}}(t,\cdot)_{\#}\tilde{\zeta}_\lambda(t,\cdot)\mathscr{L}^2, \tag{4.26}$$

and observe that, by Section 2.3, for every $q \in \mathbb{N}$, we have

$$\begin{aligned}&\sup_{t\in[0,2]}\Big(\|\rho^q_{\lambda,\boldsymbol{w}}(t,\cdot)\|_{L^\infty_x} + \|\tilde{\rho}^q_{\lambda,\boldsymbol{w}}(t,\cdot)\|_{L^\infty_x} + \|\rho_{\lambda,\boldsymbol{w}}(t,\cdot)\|_{L^\infty_x} + \|\tilde{\rho}_{\lambda,\boldsymbol{w}}(t,\cdot)\|_{L^\infty_x}\Big) \\ &\leq 4\exp\Big[\int_0^2 \|[\operatorname{div}_x \boldsymbol{w}(s,\cdot)]\|_{L^\infty_x} ds\Big].\end{aligned} \tag{4.27}$$

We finish this section by the following lemma.

**Lemma 4.8.** *$\rho^q_{\lambda,\boldsymbol{w}}$ and $\tilde{\rho}^q_{\lambda,\boldsymbol{w}}$ converge in $C([0,2]; w^* - L^\infty(\mathbb{R}^2))$ to $\rho_{\lambda,\boldsymbol{w}}$ and $\tilde{\rho}_{\lambda,\boldsymbol{w}}$ respectively as $q \to +\infty$.*

**Proof.** Let $\phi \in L^1(\mathbb{R}^2)$ and $t \in [0,2]$. Then, we have

$$\begin{aligned}\lim_{q\to+\infty}\int_{\mathbb{R}^2}\rho^q_{\lambda,\boldsymbol{w}}(t,x)\phi(x)dx &= \lim_{q\to+\infty}\int_{\mathbb{R}^2}\phi(\boldsymbol{X}_{\boldsymbol{w}}(t,x))\zeta^q_\lambda(t,x)dx \\ &= \int_{\mathbb{R}^2}\phi(\boldsymbol{X}_{\boldsymbol{w}}(t,x))\zeta_\lambda(t,x)dx \\ &= \int_{\mathbb{R}^2}\phi(x)\rho_{\lambda,\boldsymbol{w}}(t,x)dx,\end{aligned} \tag{4.28}$$

where we have used Lemma 3.1 in the second to last equality. As $\phi$ was arbitrary in $L^1(\mathbb{R}^2)$ and $t$ was arbitrary in $[0,2]$, we have shown that $\rho^q_{\lambda,\boldsymbol{w}}$ converges in $C([0,2]; w^* - L^\infty(\mathbb{R}^2))$ to $\rho_{\lambda,\boldsymbol{w}}$. A similar argument shows that $\tilde{\rho}^q_{\lambda,\boldsymbol{w}}$ converges in $C([0,2]; w^* - L^\infty(\mathbb{R}^2))$ to $\tilde{\rho}_{\lambda,\boldsymbol{w}}$. $\square$

Since $\boldsymbol{b}^q_{\lambda,\boldsymbol{w}}$ and $\tilde{\boldsymbol{b}}^q_{\lambda,\boldsymbol{w}}$ converge in $L^1_{loc}$ to $\boldsymbol{b}_{\lambda,\boldsymbol{w}}$, we may pass into the limit in Definition 2.1. Hence $\rho_{\lambda,\boldsymbol{w}}$ and $\tilde{\rho}_{\lambda,\boldsymbol{w}}$ are bounded weak solutions of (BVP) along $\boldsymbol{b}_{\lambda,\boldsymbol{w}}$ with boundary value $1/2$ at time 1, as well as bounded weak solutions of (IVP) along $\boldsymbol{b}_{\lambda,\boldsymbol{w}}$ with initial datum $\bar{\rho}$. Notice that $\rho_{\lambda,\boldsymbol{w}}$ and $\tilde{\rho}_{\lambda,\boldsymbol{w}}$ are obviously distinct.

## 5. Regularised vector fields

In this section, we conclude the proof Theorem 1.3. We will give two regularisations of the vector fields $\boldsymbol{b}^q_{\lambda,\boldsymbol{w}}$ along which the unique solutions will converge to $\rho_{\lambda,\boldsymbol{w}}$ and $\tilde{\rho}_{\lambda,\boldsymbol{w}}$ respectively. Let $\theta \in C^\infty_c(\mathbb{R}^2)$ and $\eta \in C^\infty_c(\mathbb{R})$ be standard mollifiers. We set $\theta^k(x) = k^2\theta(kx)$ and $\eta^k(t) = k\eta(kt)$. Given a vector field $\boldsymbol{w}$ in $L^\infty((0,2); C^2_c(\mathbb{R}^2;\mathbb{R}^2))$, we set $\boldsymbol{w}^k = \boldsymbol{w} \star_t \boldsymbol{\eta}^k$, which is then a bounded vector field in $C^2([0,2]\times\mathbb{R}^2;\mathbb{R}^2)$ and we denote by $\boldsymbol{X}_{\boldsymbol{w}^k}$ the flow along $\boldsymbol{w}^k$ starting at time 1. We now define

$$\boldsymbol{b}^{q,k}_{\lambda,\boldsymbol{w}}(t,x) := D_x\boldsymbol{X}_{\boldsymbol{w}^k}(t, \boldsymbol{X}^{-1}_{\boldsymbol{w}^k}(t,x)).(\boldsymbol{b}^q_\lambda \star_x \theta^k)(t, \boldsymbol{X}^{-1}_{\boldsymbol{w}^k}(t,x)) + \boldsymbol{w}^k(t,x),$$

and also

$$\tilde{\boldsymbol{b}}^{q,k}_{\lambda,\boldsymbol{w}}(t,x) := D_x\boldsymbol{X}_{\boldsymbol{w}^k}(t, \boldsymbol{X}^{-1}_{\boldsymbol{w}^k}(t,x)).(\tilde{\boldsymbol{b}}^q_\lambda \star_x \theta^k)(t, \boldsymbol{X}^{-1}_{\boldsymbol{w}^k}(t,x)) + \boldsymbol{w}^k(t,x).$$

As the above defined vector fields are smooth, they admit a unique flow: we denote by $\boldsymbol{Y}^k_s$ the flow along $\boldsymbol{b}^{q,k}_{\lambda,\boldsymbol{w}}$ starting at time $s$, and $\boldsymbol{Z}^k_s$ the flow along $\boldsymbol{b}^q_\lambda \star_x \theta^k$ starting at time $s$. By the chain rule, we have for every $t \in [0,2]$ and every $x \in \mathbb{R}^2$ that

$$\boldsymbol{Y}^k_1(t,x) = \boldsymbol{X}_{\boldsymbol{w}^k}(t, \boldsymbol{Z}^k_1(t,x)).$$

As $\boldsymbol{b}^q_\lambda \star_x \theta^k$ is divergence-free, we have

$$\boldsymbol{Z}^k_1(t,\cdot)_\# \mathscr{L}^d = \mathscr{L}^d.$$

Therefore for every $t \in [0,2]$ it holds that

$$\exp\Big(-\int_0^2 \|\operatorname{div}_x \boldsymbol{w}(s,\cdot)\|_{L^\infty_x} ds\Big)\mathscr{L}^d \leq \boldsymbol{Y}^k_1(t,\cdot)_\#\mathscr{L}^d \leq \exp\Big(\int_0^2 \|\operatorname{div}_x \boldsymbol{w}(s,\cdot)\|_{L^\infty_x} ds\Big)\mathscr{L}^d. \tag{5.1}$$

Since $\boldsymbol{Y}^k_0(t,x) = \boldsymbol{Y}^k_1(t, (\boldsymbol{Y}^k_1)^{-1}(0,x))$, we have for every $t \in [0,2]$ that

$$\exp\Big(-3\int_0^2 \|\operatorname{div}_x \boldsymbol{w}(s,\cdot)\|_{L^\infty_x} ds\Big)\mathscr{L}^d \le \boldsymbol{Y}^k_0(t,\cdot)_\#\mathscr{L}^d$$
$$\le \exp\Big(3\int_0^2 \|\operatorname{div}_x \boldsymbol{w}(s,\cdot)\|_{L^\infty_x} ds\Big)\mathscr{L}^d. \tag{5.2}$$

Let us now conclude with the proof of the main theorem.

**Proof of Theorem 1.3.** By Lemma 4.2, $D$ defined by (4.2) is dense and by Remark 4.1, it is not bounded at any time in $[0,2]$. Let $\boldsymbol{b}\in D$. By definition of the set $D$, for some $\lambda\in\mathbb{N}$ and $\boldsymbol{w}\in L^\infty((0,2);C^2_c(\mathbb{R}^2;\mathbb{R}^2))$, we have $\boldsymbol{b}=\boldsymbol{b}_{\lambda,\boldsymbol{w}}$, where $\boldsymbol{b}_{\lambda,\boldsymbol{w}}$ is defined in (4.1). We also set $\bar\rho$ according to Section 4.3.

For every $q\in\mathbb{N}$, we denote by $\rho^{q,k}_{\lambda,\boldsymbol{w}}$ and $\tilde\rho^{q,k}_{\lambda,\boldsymbol{w}}$ the unique bounded weak solution of (IVP) along $\boldsymbol{b}^{q,k}_{\lambda,\boldsymbol{w}}$ and $\tilde{\boldsymbol{b}}^{q,k}_{\lambda,\boldsymbol{w}}$ respectively, with initial datum $\bar\rho$. Recall also that $\rho^q_{\lambda,\boldsymbol{w}}$ and $\tilde\rho^q_{\lambda,\boldsymbol{w}}$ are the bounded weak solutions of (BVP) along $\boldsymbol{b}^q_{\lambda,\boldsymbol{w}}$ and $\tilde{\boldsymbol{b}}^q_{\lambda,\boldsymbol{w}}$ with boundary datum $\zeta^q(1-2^{-q},\cdot)$ and $\zeta^q(1-2^{-q-2},\cdot)$ respectively as described in Section 4.3.

**Step 1** (Choosing two regularisations $(\boldsymbol{b}^q)_{q\in\mathbb{N}}$ and $(\tilde{\boldsymbol{b}}^q)_{q\in\mathbb{N}}$ of $\boldsymbol{b}$) Let $q\in\mathbb{N}$. Note that $\boldsymbol{b}^q_{\lambda,\boldsymbol{w}}$ and $\tilde{\boldsymbol{b}}^q_{\lambda,\boldsymbol{w}}$ are both bounded, in $L^1((0,2);BV_{loc}(\mathbb{R}^2;\mathbb{R}^2))$, and nearly incompressible. Also $(\boldsymbol{b}^{q,k}_{\lambda,\boldsymbol{w}})_{k\in\mathbb{N}}$ and $(\tilde{\boldsymbol{b}}^{q,k}_{\lambda,\boldsymbol{w}})_{k\in\mathbb{N}}$ are regularisations of $\boldsymbol{b}^q_{\lambda,\boldsymbol{w}}$ and $\tilde{\boldsymbol{b}}^q_{\lambda,\boldsymbol{w}}$ respectively, which satisfy the hypothesis of Theorem 2.3 by (5.2). So by Theorem 2.3, there exists $k_q\in\mathbb{N}$ sufficiently large that

$$\begin{aligned}
&\sup_{t\in[0,2]}\int_{B_{2^q}(0)} |\rho^{q,k_q}_{\lambda,\boldsymbol{w}}(t,x)-\rho^q_{\lambda,\boldsymbol{w}}(t,x)|dx < 2^{-q},\\
&\sup_{t\in[0,2]}\int_{B_{2^q}(0)} |\tilde\rho^{q,k_q}_{\lambda,\boldsymbol{w}}(t,x)-\tilde\rho^q_{\lambda,\boldsymbol{w}}(t,x)|dx < 2^{-q}.
\end{aligned} \tag{5.3}$$

We set $\boldsymbol{b}^q := \boldsymbol{b}^{q,k_q}_{\lambda,\boldsymbol{w}}$ and $\tilde{\boldsymbol{b}}^q := \tilde{\boldsymbol{b}}^{q,k_q}_{\lambda,\boldsymbol{w}}$.

**Step 2** (Convergence of $\rho^{q,k_q}_{\lambda,\boldsymbol{w}}$ to $\rho_{\lambda,\boldsymbol{w}}$) We want to show that

$$\lim_{q\to+\infty}\sup_{t\in[0,2]}\Big|\int_{\mathbb{R}^2}(\rho^{q,k_q}_{\lambda,\boldsymbol{w}}(t,x)-\rho_{\lambda,\boldsymbol{w}}(t,x))\phi(x)dx\Big| = 0 \qquad \forall\phi\in L^1(\mathbb{R}^2). \tag{5.4}$$

Using Hölder inequality, we have

$$\Big|\int_{\mathbb{R}^2}(\rho^{q,k_q}_{\lambda,\boldsymbol{w}}(t,x)-\rho_{\lambda,\boldsymbol{w}}(t,x))(\phi(x)-\psi(x))dx\Big| \le (\|\rho^{q,k_q}_{\lambda,\boldsymbol{w}}\|_{L^\infty_x}+\|\rho_{\lambda,\boldsymbol{w}}\|_{L^\infty_x})\|\phi-\psi\|_{L^1_x}, \tag{5.5}$$

and by (5.2), the compressibility constant of the unique flow along $\boldsymbol{b}^q$ is bounded by

$$\exp\Big(\int_0^2 \|\operatorname{div}_x \boldsymbol{w}(s,\cdot)\|_{L^\infty_x} ds\Big).$$

Therefore along with equation (4.27), this yields

$$\sup_{q\in\mathbb{N}} \|\rho^{q,k_q}_{\lambda,\boldsymbol{w}}\|_{L^\infty_x} + \|\rho_{\lambda,\boldsymbol{w}}\|_{L^\infty_x} < +\infty,$$

and since $C_c(\mathbb{R}^2)$ is dense in $L^1(\mathbb{R}^2)$, in view of (5.5), to prove (5.4), we only have to prove

$$\lim_{q\to+\infty} \sup_{t\in[0,2]} \Big|\int_{\mathbb{R}^2} (\rho^{q,k_q}_{\lambda,\boldsymbol{w}}(t,x) - \rho_{\lambda,\boldsymbol{w}}(t,x))\phi(x)dx\Big| = 0 \qquad \forall \phi \in C_c(\mathbb{R}^2). \tag{5.6}$$

Accordingly, let $\phi \in C_c(\mathbb{R}^2)$, $\varepsilon > 0$, and $t \in [0,2]$. By Lemma 4.8, we can choose $\tilde{q} \in \mathbb{N}$ sufficiently large, uniformly in $t \in [0,2]$, such that for every $q \geq \tilde{q}$, we have

$$\Big|\int_{\mathbb{R}^2} (\rho^{q}_{\lambda,\boldsymbol{w}}(t,x) - \rho_{\lambda,\boldsymbol{w}}(t,x))\phi(x)dx\Big| < \varepsilon, \tag{5.7}$$

and also that

$$2^{-\tilde{q}} < \varepsilon \qquad \text{and} \qquad \operatorname{supp}\phi \subset B_{2^{\tilde{q}}}(0). \tag{5.8}$$

We then have that for every $q \geq \tilde{q}$

$$\begin{aligned} &\Big|\int_{\mathbb{R}^2} \big(\rho^{q,k_q}_{\lambda,\boldsymbol{w}}(t,x) - \rho_{\lambda,\boldsymbol{w}}(t,x)\big)\phi(x)dx\Big| \\ &\leq \Big|\int_{\mathbb{R}^2} (\rho^{q,k_q}_{\lambda,\boldsymbol{w}}(t,x)dx - \rho^{q}_{\lambda,\boldsymbol{w}}(t,x))\phi(x)dx\Big| + \Big|\int_{\mathbb{R}^2} (\rho^{q}_{\lambda,\boldsymbol{w}}(t,x) - \rho_{\lambda,\boldsymbol{w}}(t,x))\phi(x)dx\Big| \\ &\leq 2^{-q}\|\phi\|_{L^\infty_x} + \varepsilon < (\|\phi\|_{L^\infty_x} + 1)\varepsilon, \end{aligned} \tag{5.9}$$

where in the second to last inequality, we have performed a Hölder inequality on the first term we have used (5.3) and (5.8), and on the second term we have used (5.7). As $\varepsilon$ was arbitrary (5.6) follows. As $\phi \in C_c(\mathbb{R}^2)$ was arbitrary, it follows that $\rho^{q,k_q}_{\lambda,\boldsymbol{w}}$ converges to $\rho_{\lambda,\boldsymbol{w}}$ in $C([0,2]; w^* - L^\infty(\mathbb{R}^2))$ as $q \to +\infty$.

An entirely similar analysis to Step 2 also shows that $\tilde{\rho}^{q,k_q}_{\lambda,\boldsymbol{w}}$ converges to $\tilde{\rho}_{\lambda,\boldsymbol{w}}$ in $C([0,2]; w^* - L^\infty(\mathbb{R}^2))$ as $q \to +\infty$. As $\rho^{q,k_q}_{\lambda,\boldsymbol{w}}$ and $\tilde{\rho}^{q,k_q}_{\lambda,\boldsymbol{w}}$ are the unique bounded weak solutions of (IVP) along $\boldsymbol{b}^q$ and $\tilde{\boldsymbol{b}}^q$ respectively with initial datum $\bar{\rho}$, the thesis is proven. □

## Data availability

No data was used for the research described in the article.